# NORMALIZED LEAST-SQUARES ESTIMATION IN TIME-VARYING ARCH MODELS

By Piotr Fryzlewicz,[1] Theofanis Sapatinas
and Suhasini Subba Rao[2]

*University of Bristol, University of Cyprus
and Texas A&M University*

We investigate the time-varying ARCH (tvARCH) process. It is shown that it can be used to describe the slow decay of the sample autocorrelations of the squared returns often observed in financial time series, which warrants the further study of parameter estimation methods for the model.

Since the parameters are changing over time, a successful estimator needs to perform well for small samples. We propose a kernel normalized-least-squares (kernel-NLS) estimator which has a closed form, and thus outperforms the previously proposed kernel quasi-maximum likelihood (kernel-QML) estimator for small samples. The kernel-NLS estimator is simple, works under mild moment assumptions and avoids some of the parameter space restrictions imposed by the kernel-QML estimator. Theoretical evidence shows that the kernel-NLS estimator has the same rate of convergence as the kernel-QML estimator. Due to the kernel-NLS estimator's ease of computation, computationally intensive procedures can be used. A prediction-based cross-validation method is proposed for selecting the bandwidth of the kernel-NLS estimator. Also, we use a residual-based bootstrap scheme to bootstrap the tvARCH process. The bootstrap sample is used to obtain pointwise confidence intervals for the kernel-NLS estimator. It is shown that distributions of the estimator using the bootstrap and the "true" tvARCH estimator asymptotically coincide.

We illustrate our estimation method on a variety of currency exchange and stock index data for which we obtain both good fits to the data and accurate forecasts.

Received June 2007; revised June 2007.

[1]Supported in part by the Department of Mathematics, Imperial College, UK.

[2]Supported by the DFG (DA 187/12-2).

**1. Introduction.** Among models for log-returns $X_t = \log(P_t/P_{t-1})$ on speculative prices $P_t$ (such as currency exchange rates, share prices, stock indices, etc.), the stationary ARCH($p$) [Engle (1982)] and GARCH($p$, $q$) [Bollerslev (1986) and Taylor (1986)] processes have gained particular popularity and have become standard in the financial econometrics literature as they model well the volatility of financial markets over short periods of time. For a review of recent advances on those and related models, we refer the reader to Fan and Yao (2003) and Giraitis, Leipus and Surgailis (2005).

The modeling of financial data using nonstationary time series models has recently attracted considerable attention. Arguments for using such models were laid out, for example, in Fan, Jiang, Zhang and Zhou (2003), Mikosch and Stărică (2000, 2003, 2004), Mercurio and Spokoiny (2004a, 2004b), Stărică and Granger (2005) and Fryzlewicz et al. (2006).

Recently, Dahlhaus and Subba Rao (2006) generalized the class of ARCH($p$) processes to include processes whose parameters were allowed to change "slowly" through time. The resulting model, called the time-varying ARCH($p$) [tvARCH($p$)] process, is defined as

$$(1) \qquad X_{t,N} = \sigma_{t,N} Z_t, \qquad \sigma_{t,N}^2 = a_0\left(\frac{t}{N}\right) + \sum_{j=1}^{p} a_j\left(\frac{t}{N}\right) X_{t-j,N}^2,$$

for $t = 1, 2, \ldots, N$, where $\{Z_t\}_t$ are independent and identically distributed random variables with $\mathbb{E}(Z_t) = 0$ and $\mathbb{E}(Z_t^2) = 1$. In this paper, we focus on how the tvARCH($p$) process can be used to characterize some of the features present in financial data, estimation methods for small samples, bootstrapping the tvARCH($p$) process and the fitting of the tvARCH($p$) process to data.

In Section 2, we show how the tvARCH($p$) process can be used to describe the slow decay of the sample autocorrelations of the squared returns often observed in financial log-returns and usually attributed to the long memory of the underlying process. This is despite the true nonstationary correlations decaying geometrically fast to zero. Thus, the tvARCH($p$) process, due to its nonstationarity, captures the appearance of long memory which is present in many financial datasets: a feature also exhibited by a short memory GARCH($1, 1$) process with structural breaks [Mikosch and Stărică (2000, 2003, 2004)—note that this effect goes back to Bhattacharya, Gupta and Waymire (1983)].

The benchmark method for the estimation of stationary ARCH($p$) parameters is the quasi-maximum likelihood (QML) estimator. Motivated by this, Dahlhaus and Subba Rao (2006) use a localized kernel-based quasi-maximum likelihood (kernel-QML) method for estimating the parameters of a tvARCH($p$) process. However, the kernel-QML estimator for small sample sizes is not very reliable, since the QML tends to be shallow about the minimum for small sample sizes [Shephard ([1996](#)) and Bose and Mukherjee ([2003](#))]. This is of particular relevance to tvARCH($p$) processes, where in regions of nonstationarity, we need to base our estimator on only a few observations to avoid a large bias. Furthermore, the parameter space of the estimator is restricted to $\inf_j a_j(u) > 0$. However, it is suggested in the examples in Section [6](#) that over large periods of time some of the higher-order parameters should be zero. This renders the assumption $\inf_j a_j(u) > 0$ rather unrealistic. In addition, evaluation of the kernel-QML estimator at every time point is computationally quite intensive. Therefore, bandwidth selection based on a data driven procedure, where the kernel-QML estimator has to be evaluated at each time point for different bandwidths, may not be feasible for even moderately large sample sizes.

A rival class of estimators are least-squares-based and are known to have good small-sample properties [Bose and Mukherjee ([2003](#))]. These types of estimators will be the focal point in this paper. In Section [3](#) and the following sections, we propose and thoroughly analyze a (suitably localized and normalized) least-squares-type estimator for the tvARCH($p$) process which, unlike the kernel-QML estimator mentioned above, enjoys the following properties: (i) very good performance for small samples, (ii) simplicity and closed form and (iii) rapid computability. In addition, it does allow $\inf_j a_j(u) = 0$, thereby avoiding the parameter space restriction described above.

In Section [3.1](#), we consider a general class of localized weighted least-squares estimators for tvARCH($p$) process and study their sampling properties. We show that their small sample performance, sampling properties and moment assumptions depend on the weight function used.

In Section [3.3](#), we investigate weight functions that lead to estimators which are close to the kernel-QML estimator for large samples and easy to compute. In fact, we show that the weight functions which have the most desirable properties contain unknown parameters. This motivates us in Section [3.4](#) to propose the two-stage kernel normalized-least-squares (kernel-NLS) estimator where in the first stage we estimate the weight function which we use in the second stage as the weight in the least-squares estimator. The two-stage kernel-NLS estimator has the same sampling properties as if the true weight function were a priori known, and has the same rate of convergence as the kernel-QML estimator. In Section [3.6](#), we state some of the results from extensive simulation studies which show that for small sample sizes the two-stage kernel-NLS estimator performs better than the kernel-QML estimator. This suggests that at least in the nonstationary setup, the two-stage kernel-NLS estimator is a viable alternative to the kernel-QML estimator.

In Section [4](#), we propose a cross-validation method for selecting the bandwidth of the two-stage kernel-NLS estimator. The proposed cross-validation procedure for tvARCH($p$) processes is based on one-step-ahead prediction of the data to select the bandwidth. The closed form solution of the two-stage kernel-NLS estimator means that, for every bandwidth, the estimator can be evaluated rapidly. The computation ease of the two-stage kernel-NLS estimator means that it is simple to implement a cross-validation method based on this scheme. We discuss some of the implementation issues associated with the procedure and show that its computational complexity remains low.

In Section 5, we bootstrap the tvARCH($p$) process. This allows us to obtain finite sample pointwise confidence intervals for the tvARCH($p$) parameter estimators. The scheme is based on bootstrapping the estimated residuals, which we use, together with the estimated tvARCH($p$) parameters, to construct the bootstrap sample. Again, the fact that the bootstrapping scheme is computationally feasible is only due to the rapid computability of the two-stage kernel-NLS estimator. We show that the distribution of the bootstrap tvARCH($p$) estimator asymptotically coincides with the "true" tvARCH($p$) estimator. The method and results in this section may also be of independent interest.

In Section 6, we demonstrate that our estimation methodology gives a very good fit to data for the USD/GBP currency exchange and FTSE stock index datasets, and we also exhibit bootstrap pointwise confidence intervals for the estimated parameters. In Section 7, we test the long-term volatility forecasting ability of the tvARCH($p$) process with $p = 0, 1, 2$, where the parameters are estimated via the two-stage kernel-NLS estimator. We show that, for a variety of currency exchange datasets, our forecasting methodology outperforms the stationary GARCH(1, 1) and EGARCH(1, 1) techniques. However, it is interesting to observe that the latter two methods give slightly superior results for a selection of stock index datasets.

Proofs of the results in the paper are outlined in the Appendix. Further details of the proofs can be found in the accompanying technical report, available from the authors or from http://www.maths.bris.ac.uk/ mapzf/ tvarch/trNLS.pdf.

**2. The tvARCH($p$) process: preliminary results and motivation.** In this section, we discuss some of the properties of the tvARCH($p$) process.

2.1. *Notation, assumptions and main ingredients.* We first state the assumptions used throughout the paper.

ASSUMPTION 1. Suppose $\{X_{t,N}\}_t$ is a tvARCH ($p$) process. We assume that the time-varying parameters $\{a_j(u)\}_j$ and the innovations $\{Z_t\}_t$ satisfy the following conditions:

(i) There exist $0 < \rho_1 \leq \rho_2 < \infty$ and $0 < \delta < 1$ such that, for all $u \in (0, 1]$, $\rho_1 \leq a_0(u) \leq \rho_2$, and $\sup_u \sum_{j=1}^p a_j(u) \leq 1 - \delta$.

(ii) There exist $\beta \in (0, 1]$ and a finite constant $K > 0$ such that for $u, v \in (0, 1]$

$$|a_j(u) - a_j(v)| \leq K|u - v|^\beta \qquad \text{for each} \quad j = 0, 1, \ldots, p.$$

(iii) For some $\gamma > 0$, $\mathbb{E}(|Z_t|^{4(1+\gamma)}) < \infty$.

(iv) For some $\eta > 0$ and $0 < \delta < 1$, $m_{1+\eta} \sup_u \sum_{j=1}^p a_j(u) \leq 1 - \delta$, where $m_{1+\eta} = \{\mathbb{E}(|Z_t|^{2(1+\eta)})\}^{1/(1+\eta)}$.

Assumption 1(i) implies that $\sup_{t,N} \mathbb{E}(X_{t,N}^2) < \infty$. Assumption 1(i), (ii) means that the tvARCH($p$) process can locally be approximated by a stationary process. We require Assumption 1(iii), (iv) to show asymptotic normality of the two-stage kernel-NLS estimator (defined in Section 3.4). Comparing $m_{1+\eta} \sup_u \sum_{j=1}^p a_j(u) \leq 1 - \delta$ with the assumption required to show asymptotic normality of the kernel-QML estimator ($m_1 \sup_u \sum_{j=1}^p a_j(u) \leq 1 - \delta$, where we note that $m_1 = 1$), it is only a mildly stronger assumption, as we only require it to hold for *some* $\eta > 0$. In other words, if the moment function $m_\nu$ increases smoothly with $\nu$, and $m_1 \sup_u \sum_{j=1}^p a_j(u) \leq 1 - \delta$, then there exists a $\eta > 0$ and $0 < \delta_1 < 1$ such that $m_{1+\eta} \sup_u \sum_{j=1}^p a_j(u) \leq 1 - \delta_1$ [which satisfies Assumption 1(iv)].

In order to prove results concerning the tvARCH($p$) process, Dahlhaus and Subba Rao (2006) define the stationary process $\{\tilde{X}_t(u)\}_t$. Let $u \in (0, 1]$ and suppose that, for each fixed $u$, $\{\tilde{X}_t(u)\}_t$ satisfies the model

$$(2) \qquad \tilde{X}_t(u) = \tilde{\sigma}_t(u) Z_t, \qquad \tilde{\sigma}_t^2(u) = a_0(u) + \sum_{j=1}^p a_j(u) \tilde{X}_{t-j}^2(u).$$

The following lemma is a special case of Corollary 4.2 in Subba Rao (2006), where it was shown that $\{\tilde{X}_t^2(u)\}_t$ can be regarded as a stationary approximation of the nonstationary process $\{X_{t,N}^2\}_t$ about $u \approx t/N$, which is why $\{X_{t,N}\}_t$ can be regarded as a *locally stationary process*. We can treat the lemma below as the stochastic version of Hölder continuity.

LEMMA 1. *Suppose $\{X_{t,N}\}_t$ is a tvARCH($p$) process which satisfies Assumption 1(i), (ii), and let $\{\tilde{X}_t(u)\}_t$ be defined as in (2). Then, for each fixed $u \in (0, 1]$, we have that $\{\tilde{X}_t^2(u)\}_t$ is a stationary, ergodic process such that*

$$(3) \qquad |X_{t,N}^2 - \tilde{X}_t^2(u)| \leq \frac{1}{N^\beta} V_{t,N} + \left|u - \frac{t}{N}\right|^\beta W_t \qquad \text{almost surely,}$$

*and $|\tilde{X}_t^2(u) - \tilde{X}_t^2(v)| \leq |u - v|^\beta W_t$, almost surely, where $\{V_{t,N}\}_t$ and $\{W_t\}_t$ are well-defined positive processes, and $\{W_t\}_t$ is a stationary process. In addition, if we assume that Assumption 1(iv) holds, then we have $\sup_{t,N} \mathbb{E}|V_{t,N}|^{1+\eta} < \infty$ and $\mathbb{E}|W_t|^{1+\eta} < \infty$.*

Several of the estimators considered in this paper [e.g., the estimators defined in (4) and (7), etc.] are local or global averages of functions of the tvARCH($p$) process. Unlike stationary ARCH($p$) (or more general stationary) processes, we cannot study the sampling properties of these estimators by simply letting the sample size grow. Instead, we use the rescaling by $N$ to obtain a meaningful asymptotic theory. The underlying principle to studying an estimator at a particular time $t$, is to keep the ratio $t/N$ fixed and let $N \to \infty$ [Dahlhaus (1997)]. However, the tvARCH($p$) process varies for different $N$, which is the reason for introducing the stationary approximation.

Throughout the paper, $\xrightarrow{\mathcal{P}}$ and $\xrightarrow{\mathcal{D}}$ denote convergence in probability and in distribution, respectively.

2.2. *The covariance structure and the long memory effect.* The following proposition shows the behavior of the true autocovariance function of the squares of a tvARCH($p$) process.

PROPOSITION 1. *Suppose $\{X_{t,N}\}_t$ is a* tvARCH($p$) *process which satisfies Assumption* 1(i), (ii), *and assume that $\{\mathbb{E}(Z_t^4)\}^{1/2} \sup_u \sum_{j=1}^p a_j(u) \le 1 - \delta$, for some $0 < \delta < 1$. Then, for some $\rho \in (1-\delta, 1)$ and a fixed $h \ge 0$, we have*

$$\sup_{t,N} |\operatorname{cov}\{X_{t,N}^2, X_{t+h,N}^2\}| \le K\rho^h,$$

*for some finite constant $K > 0$ that is independent of $h$.*

If the fourth moment of the process $\{X_{t,N}\}_t$ exists, then Proposition 1 implies that $\{X_{t,N}^2\}_t$ is a short memory process.

However, we now show that the sample autocovariance of the process $\{X_{t,N}^2\}_t$, computed under the wrong premise of stationarity, does not necessarily decay to zero. Typically, if we believed that the process $\{X_{t,N}^2\}_t$ were stationary, we would use $\mathcal{S}_N(h)$ as an estimator of $\operatorname{cov}\{X_{t,N}^2, X_{t+h,N}^2\}$, where

$$(4) \qquad \mathcal{S}_N(h) = \frac{1}{N-h} \sum_{t=1}^{N-h} X_{t,N}^2 X_{t+h,N}^2 - (\bar{X}_N)^2$$

and

$$\bar{X}_N = \frac{1}{N-h} \sum_{t=1}^{N-h} X_{t,N}^2.$$

Denote $\mu(u) = \mathbb{E}(\tilde{X}_t^2(u))$ and $c(u,h) = \operatorname{cov}\{\tilde{X}_t^2(u), \tilde{X}_{t+h}^2(u)\}$ for each $u \in (0,1]$ and $h \ge 0$.

The following proposition shows the behavior of the sample autocovariance of the squares of a tvARCH($p$) process, evaluated under the wrong assumption of stationarity.

PROPOSITION 2. *Suppose $\{X_{t,N}\}_t$ is a tvARCH($p$) process which satisfies Assumption 1(i), (ii), and assume that, for some $0 < \zeta \leq 2$ and $0 < \delta < 1$, $\{\mathbb{E}(|Z_t|^{2(2+\zeta)})\}^{1/(2+\zeta)} \sup_u \sum_{j=1}^p a_j(u) \leq 1 - \delta$. Then, for fixed $h > 0$, as $N \to \infty$, we have*

$$
(5) \qquad \mathcal{S}_N(h) \xrightarrow{\mathcal{P}} \int_0^1 c(u,h)\, du + \int \int_{\{0 \leq u < v \leq 1\}} \{\mu(u) - \mu(v)\}^2 \, du \, dv.
$$

According to Proposition 2, since the autocovariance of the squares of a tvARCH($p$) process decays to zero exponentially fast as $h \to \infty$, so does the first integral in (5). However, the appearance of persistent correlations would still appear if the second integral were nonzero. We consider the simple example when the mean of the squares increases linearly, that is, if $\mu(u) = cu$, for some nonzero constant $c$. In this case, the second integral in (5) reduces to $c^2/12$. In other words, the long memory effect is due to changes in the unconditional variance of the tvARCH($p$) process.

## 3. The kernel-NLS estimator and its asymptotic properties.

Typically, to estimate the parameters of a stationary ARCH($p$) process, a QML estimator is used, where the likelihood is constructed as if the innovations were Gaussian. The main advantage of the QML estimator is that, even in the case that the innovations are non-Gaussian, it is consistent and asymptotically normal. In contrast, Straumann (2005) has shown that under misspecification of the innovation distribution, the resulting non-Gaussian maximum likelihood estimator is inconsistent. As it is almost impossible to specify the distribution of the innovations, this makes the QML estimator the benchmark method when estimating stationary ARCH($p$) parameters.

A localized version of the QML estimator is used to estimate the parameters of a tvARCH($p$) process in Dahlhaus and Subba Rao (2006). To prove the sampling results, the asymptotics are done in the rescaled time framework. In practice, a good estimator is obtained if the process is close to stationary over a relatively large region. However, the story is completely different over much shorter regions. As noted in the Section 1, in estimation over a short period of time (which will often be the case for nonstationary processes), the performance of the QML estimator is quite poor.

Rival methods are least-squares-type estimators which are known to have good small sample properties. In this section, we focus on kernel weighted least-squares as a method for estimating the parameters of a tvARCH($p$) process. To this end, we define the kernel $W : [-1/2, 1/2] \to \mathbb{R}$, which is a function of bounded variation and satisfies the standard conditions: $\int_{-1/2}^{1/2} W(x)\, dx = 1$ and $\int_{-1/2}^{1/2} W^2(x)\, dx < \infty$.

3.1. *Kernel weighted least-squares for* tvARCH $(p)$ *processes.* It is straightforward to show that the squares of the tvARCH$(p)$ process satisfy the autoregressive representation $X_{t,N}^2 = a_0(\frac{t}{N}) + \sum_{j=1}^p a_j(\frac{t}{N})X_{t-j,N}^2 + (Z_t^2 - 1)\sigma_{t,N}^2$. For reasons that will become obvious later, we weight the least squares representation with the weight function $\kappa(u_0, \mathcal{X}_{k-1,N})$, where $\mathcal{X}_{k-1,N}^T = (1, X_{k-1,N}^2, \ldots, X_{k-p,N}^2)$, and define the following weighted least-squares criterion:

$$(6) \qquad \mathcal{L}_{t_0,N}(\underline{\alpha}) = \sum_{k=p+1}^N \frac{1}{bN} W\left(\frac{t_0-k}{bN}\right) \frac{(X_{k,N}^2 - \alpha_0 - \sum_{j=1}^p \alpha_j X_{k-j,N}^2)^2}{\kappa(u_0, \mathcal{X}_{k-1,N})^2}.$$

If $|u_0 - t_0/N| < 1/N$, we use $\underline{\hat{a}}_{t_0,N}$ as an estimator of $\underline{a}(u_0) = (a_0(u), a_1(u), \ldots, a_p(u))^T$, where

$$(7) \qquad \underline{\hat{a}}_{t_0,N} = \arg\min_{\underline{a}} \mathcal{L}_{t_0,N}(\underline{a}).$$

Since $\underline{\hat{a}}_{t_0,N}$ is a least-squares estimator, it has the advantage of a closed form solution, that is, $\underline{\hat{a}}_{t_0,N} = \{\mathcal{R}_{t_0,N}\}^{-1}\underline{r}_{t_0,N}$, where

$$\mathcal{R}_{t_0,N} = \sum_{k=p+1}^N \frac{1}{bN} W\left(\frac{t_0-k}{bN}\right) \frac{\mathcal{X}_{k-1,N}\mathcal{X}_{k-1,N}^T}{\kappa(u_0, \mathcal{X}_{k-1,N})^2},$$

$$\underline{r}_{t_0,N} = \sum_{k=p+1}^N \frac{1}{bN} W\left(\frac{t_0-k}{bN}\right) \frac{X_{k,N}^2 \mathcal{X}_{k-1,N}}{\kappa(u_0, \mathcal{X}_{k-1,N})^2}.$$

3.2. *Asymptotic properties of the kernel weighted least-squares estimator.* We now obtain the asymptotic sampling properties of $\underline{\hat{a}}_{t_0,N}$.

To show asymptotic normality we require the following definitions:

$$(8) \qquad \mathcal{A}_k(u) = \frac{\tilde{\mathcal{X}}_{k-1}(u)\tilde{\mathcal{X}}_{k-1}^T(u)}{\kappa(u_0, \tilde{\mathcal{X}}_{k-1}(u))^2}, \qquad \mathcal{D}_k(u) = \frac{\tilde{\sigma}_k^4(u)\tilde{\mathcal{X}}_{k-1}(u)\tilde{\mathcal{X}}_{k-1}^T(u)}{\kappa(u_0, \tilde{\mathcal{X}}_{k-1}(u))^4}$$

and

$$(9) \qquad \begin{aligned} \mathcal{B}_{t_0,N}(\underline{\alpha}) = &\sum_{k=p+1}^N \frac{1}{bN} W\left(\frac{t_0-k}{bN}\right) \\ &\times \left[ \frac{\{X_{k,N}^2 - \alpha_0 - \sum_{j=1}^p \alpha_j X_{k-j,N}^2\}^2}{\kappa(u_0, \mathcal{X}_{k-1,N})^2} \right. \\ &\left. - \frac{\{\tilde{X}_k^2(u_0) - \alpha_0 - \sum_{j=1}^p \alpha_j \tilde{X}_{k-j}^2(u_0)\}^2}{\kappa(u_0, \tilde{\mathcal{X}}_{k-1}(u))^2} \right], \end{aligned}$$

where $\tilde{\mathcal{X}}_{t-1}(u) = (1, \tilde{X}_{t-1}^2(u), \ldots, \tilde{X}_{t-p}^2(u))$. We point out that if $\{X_{t,N}\}_t$ were a stationary process then $\mathcal{B}_{t_0,N}(\underline{\alpha}) \equiv 0$.

In the following proposition we obtain consistency and asymptotic normality of $\hat{\underline{a}}_{t_0,N}$. We denote $\nabla f(u,\underline{a}) = (\frac{\partial f(u,\underline{a})}{\partial a_0}, \ldots, \frac{\partial f(u,\underline{a})}{\partial a_p})^T$, and set $\underline{x} = (1, x_1, x_2, \ldots, x_p)$ and $\underline{y} = (1, y_1, y_2, \ldots, y_p)$.

PROPOSITION 3. *Suppose $\{X_{t,N}\}_t$ is a tvARCH($p$) process which satisfies Assumption 1(i), (ii), (iii), and let $\hat{\underline{a}}_{t_0,N}$, $\mathcal{A}_t(u)$, $\mathcal{D}_t(u)$ and $\mathcal{B}_{t_0,N}(\underline{a})$ be defined as in (7), (8) and (9), respectively. We further assume that $\kappa$ is bounded away from zero and we have a type of Lipschitz condition on the weighted least-squares; that is, for all $1 \le i \le p$, $|\frac{x_i}{\kappa(u,\underline{x})} - \frac{y_i}{\kappa(u,\underline{y})}| \le K \sum_{j=1}^{p} |x_j - y_j|$, for some finite constant $K > 0$. Also, assume for all $1 \le i \le p$ that $\sup_{k,N} \mathbb{E}(\frac{X_{k-i,N}^4}{\kappa(u_0,\mathcal{X}_{k-1,N})^2}) < \infty$, and suppose $|u_0 - t_0/N| < 1/N$.*

(i) *Then we have $\hat{\underline{a}}_{t_0,N} \xrightarrow{\mathcal{P}} \underline{a}(u_0)$, with $b \to 0, bN \to \infty$ as $N \to \infty$.*

(ii) *If in addition we assume for all $1 \le i \le p$ and some $\nu > 0$ that $\sup_{k,N} \mathbb{E}(\frac{X_{k-i,N}^{8+2\nu}}{\kappa(u_0,\mathcal{X}_{k-1,N})^{4+\nu}}) < \infty$, then we have $\nabla \mathcal{B}_{t_0,N}(\underline{a}(u_0)) = O_p(b^\beta)$ and*

$$(10) \quad \begin{aligned} &\sqrt{bN}(\hat{\underline{a}}_{t_0,N} - \underline{a}(u_0)) + \tfrac{1}{2}\sqrt{bN}\mathbb{E}[\mathcal{A}_t(u_0)]^{-1}\nabla\mathcal{B}_{t_0,N}(\underline{a}(u_0)) \\ &\xrightarrow{\mathcal{D}} \mathcal{N}(0, w_2\mu_4 \mathbb{E}[\mathcal{A}_t(u_0)]^{-1}\mathbb{E}[\mathcal{D}_t(u_0)]\mathbb{E}[\mathcal{A}_t(u_0)]^{-1}), \end{aligned}$$

*with $b \to 0$, $bN \to \infty$ as $N \to \infty$, where $w_2 = \int_{-1/2}^{1/2} W^2(x)\,dx$ and $\mu_4 = \mathrm{var}(Z_t^2)$.*

At first glance the above assumptions may appear quite technical, but we note that in the case $\kappa(\cdot) \equiv 1$, they are standard in least-squares estimation. Furthermore, if the weight function $\kappa$ is bounded away from zero and Lipschitz continuous [i.e., $\sup_{\underline{x},\underline{y}} |\kappa(u,\underline{x}) - \kappa(u,\underline{y})| \le K \sum_{j=1}^{p}|x_j - y_j|$, for some finite constant $K > 0$], then it is straightforward to see that $|\frac{x_i}{\kappa(u,\underline{x})} - \frac{y_i}{\kappa(u,\underline{y})}| \le K \sum_{j=1}^{p} |x_j - y_j|$. In the following section, we will suggest a $\kappa(\cdot)$ that is ideal for tvARCH($p$) estimation and satisfies the required conditions.

3.3. *Choice of weight function $\kappa$.* By considering both theoretical and empirical evidence, we now investigate various choices of weight functions. To do this, we study Proposition 3 and consider the $\kappa$ which yields an estimator which requires only weak moment assumptions and has minimal error [see (10)]. Considering first the bias in (10), if $\sqrt{bN}b^\beta \to 0$, then the bias converges in probability to zero. Instead we focus attention on (i) the variance $\mathbb{E}[\mathcal{A}_t(u_0)]^{-1}\mathbb{E}[\mathcal{D}_t(u_0)]\mathbb{E}[\mathcal{A}_t(u_0)]^{-1}$ and (ii) derivation under low moment assumptions.

In the stationary ARCH framework, Giraitis and Robinson (2001), Bose and Mukherjee (2003), Horváth and Liese (2004) and Ling (2007) have considered the weighted least-squares estimator for different weight functions.

Giraitis and Robinson (2001) use the Whittle likelihood to estimate the parameters of a stationary $ARCH(\infty)$ process. Adapted to the nonstationary setting, the local Whittle likelihood estimator and the local weighted least-squares estimator are asymptotically equivalent when $\kappa(\cdot) \equiv 1$. Studying their assumptions, $\sup_{t,N} \mathbb{E}(X_{t,N}^4) < \infty$ and $\sup_{t,N} \mathbb{E}(X_{t,N}^{8+2\nu}) < \infty$, for some $\nu > 0$, are required to show consistency and asymptotic normality. Assuming normality of the innovations $\{Z_t\}_t$ and interpreting these conditions in terms of the coefficients of the tvARCH($p$) process, they imply that $\sup_u \sum_{j=1}^p a_j(u) < 1/\sqrt{3}$ is required for consistency and $\sup_u \sum_{j=1}^p a_j(u) < 1/\{\mathbb{E}(Z_t^{8+2\nu})\}^{1/(4+\nu)}$ for asymptotic normality. In other words, the tvARCH($p$) process should be close to a white noise process for the sampling results to be valid.

On the other hand, Bose and Mukherjee (2003) use a two-stage least-squares procedure to estimate the stationary ARCH($p$) parameters. In the first stage, they use least-squares with weight function $\kappa(\cdot) \equiv 1$ and in the second stage—a least-squares estimator with $\kappa = \hat{\sigma}_t^2$, where $\hat{\sigma}_t^2$ is an estimator of the conditional variance. An advantage of their scheme is that, asymptotically, it has the same distribution variance as the QML estimator. However, because in the first stage they use the weight $\kappa(\cdot) \equiv 1$, their method requires the same set of assumptions as in Giraitis and Robinson (2001).

To reduce the high moment restrictions, Horváth and Liese (2004) use random weights of the form $\kappa(u, \mathcal{X}_{k-1,N}) = 1 + \sum_{j=1}^p X_{k-j,N}^2$ to estimate stationary ARCH($p$) parameters, and Ling (2007) uses a similar weighting to estimate the parameters of a stationary ARMA–GARCH process. The main advantage of using this choice of weight functions is that under Assumption 1(i), (ii), (iii) the moment assumptions in Proposition 3 are satisfied.

Motivated by the discussion above, let us consider weight functions which have the form $\kappa(u, \mathcal{X}_{k-1,N}) = g(u) + \sum_{j=1}^p \rho_j(u) X_{k-j,N}^2$. We will make some comparisons with the kernel-QML estimator considered in Dahlhaus and Subba Rao (2006), who showed that the kernel-QML estimator is asymptotically normal with variance $w_2 \mu_4 \mathbb{E}[\Sigma_t(u_0)]^{-1}$, where

$$\Sigma_t(u_0) = \frac{\tilde{\mathcal{X}}_{t-1}(u_0)^T \tilde{\mathcal{X}}_{t-1}(u_0)}{\tilde{\sigma}_t^4(u_0)}. \tag{11}$$

It is worth noting that if $\{\rho_j(u)\}$ are bounded away from zero, then the conditions in Proposition 3 are fulfilled with no additional assumptions. For the purposes of this discussion only, let us assume for a moment that $\inf_j a_j(u) > 0$ (although this is not a requirement for our estimation methodology to be valid). In order to select $g(\cdot)$ and $\rho_j(\cdot)$, we first observe that if $\underline{a}(u_0)$ were known then letting $\kappa(u_0, \mathcal{X}_{k-1,N}) = a_0(u_0) + \sum_{j=1}^p a_j(u_0) X_{k-1,N}^2$ would be the ideal choice [provided $\inf_j a_j(u_0) > 0$] as the asymptotic variance of the resulting kernel weighted least-squares estimator would be the same as the kernel-QML estimator. Clearly this weight function is unknown, and for this reason we call it the "oracle" weight. Instead, we look for a closely related alternative, which is computationally simple to evaluate and avoids the requirement that $\inf_j a_j(u_0) > 0$. Let us consider a weight function $\kappa(u, \mathcal{X}_{k-1,N}) = g(u) + \sum_{j=1}^p X_{k-j,N}^2$ [which is in the spirit of the solution proposed by Horváth and Liese ([2004](#)) for stationary ARCH($p$) processes] and compare it to the oracle weight. For convenience, the estimator using the weight function $g(u) + \sum_{j=1}^p X_{k-j,N}^2$ we call the $g$-estimator, and the estimator using the oracle weight we call the oracle estimator.

Using Proposition [3](#), we see that the asymptotic distribution variance of the $g$-estimator and the oracle estimator is $w_2 \mu_4 \mathbb{E}[\mathcal{A}_t^{(g)}(u)]^{-1} \mathbb{E}[\mathcal{D}_t^{(g)}(u)] \times \mathbb{E}[\mathcal{A}_t^{(g)}(u)]^{-1}$ and $w_2 \mu_4 \mathbb{E}[\Sigma_t(u_0)]^{-1}$, respectively, where

$$
\begin{aligned}
\mathcal{A}_k^{(g)}(u) &= \frac{\tilde{\mathcal{X}}_{k-1}(u)\tilde{\mathcal{X}}_{k-1}^T(u)}{[g(u) + \sum_{j=1}^p \tilde{X}_{k-j}^2(u)]^2}, \\
\mathcal{D}_k^{(g)}(u) &= \frac{\tilde{\sigma}_k^4(u)\tilde{\mathcal{X}}_{k-1}(u)\tilde{\mathcal{X}}_{k-1}^T(u)}{[g(u) + \sum_{j=1}^p \tilde{X}_{k-j}^2(u)]^4},
\end{aligned}
\tag{12}
$$

and $\Sigma_t(u)$ is defined in ([11](#)). Let $\alpha(u) = \sum_{j=1}^p a_j(u)$, $\beta(u) = 1/\min_{j=1}^p a_j(u)$ and $|A|_{\det}$ denote the determinant of a matrix. By bounding $\mathcal{A}_t^{(g)}(u)$ and $\mathcal{D}_t^{(g)}(u)$ from both above and below we obtain

$$
\begin{aligned}
\varpi(g)^{-4}|\mathbb{E}[\Sigma_t(u)]|_{\det}^{-1} &\le |\mathbb{E}[\mathcal{A}_t^{(g)}(u)]^{-1}\mathbb{E}[\mathcal{D}_t^{(g)}(u)]\mathbb{E}[\mathcal{A}_t^{(g)}(u)]^{-1}|_{\det} \\
&\le \varpi(g)^4 |\mathbb{E}[\Sigma_t(u)]|_{\det}^{-1},
\end{aligned}
\tag{13}
$$

where

$$
\varpi(g) = \left(\frac{a_0(u) + g(u)\,\alpha(u)}{g(u)}\right)\left(\frac{g(u) + \beta(u)a_0(u)}{a_0(u)}\right).
$$

Examining ([13](#)), we have an upper and lower bound for the asymptotic distribution variance of the $g$-estimator in terms of the asymptotic variance of the oracle estimator. It is easily seen that the difference $(\varpi(g)^4 - \varpi(g)^{-4})|\mathbb{E}[\Sigma_t(u)]|_{\det}^{-1}$ and the upper bound $\varpi(g)^4|\mathbb{E}[\Sigma_t(u)]|_{\det}^{-1}$ are minimized when $g^*(u) = (a_0(u))/([\min_{1\le j\le p} a_j(u)]\sum_{j=1}^p a_j(u))$. However, $g^*(u)$ depends on unknown parameters and is highly sensitive to small values of $a_j(u)$, hence it is inappropriate as a weight function. Instead, we consider a close relative $g(u) := \mu(u) = a_0(u)/(1 - \alpha(u))$, where $\mu(u) = \mathbb{E}[\tilde{X}_t^2(u)]$. In this case, using ([13](#)), we obtain the following upper and lower bound for the asymptotic variance of the kernel weighted least-squares estimator in terms of the oracle variance:

$$
\begin{aligned}
|\mathbb{E}[\Sigma_t(u)]|_{\det}^{-1}\omega(u)^{-1} &\le |\mathbb{E}[\mathcal{A}_t^{(\mu)}(u)]^{-1}\mathbb{E}[\mathcal{D}_t^{(\mu)}(u)]\mathbb{E}[\mathcal{A}_t^{(\mu)}(u)]^{-1}|_{\det} \\
&\le |\mathbb{E}[\Sigma_t(u)]|_{\det}^{-1}\omega(u),
\end{aligned}
\tag{14}
$$

where

$$\omega(u) = \left(\frac{1 + \beta(u)[1 - \alpha(u)]}{1 - \alpha(u)}\right)^4.$$

We notice that the upper and lower bounds in (14) do not depend on the magnitude of $a_0(u)$.

Since $\frac{a_0(u)}{1-\alpha(u)} = \mathbb{E}(\tilde{X}_t^2(u)) = \mu(u)$, which is the local mean, it can easily be estimated from $\{X_{k,N}\}$. In the following section, we use it to estimate the weight function $\kappa(u_0, \mathcal{X}_{k-1,N}) = \mu(u_0) + \mathcal{S}_{k-1,N}$, where $\mathcal{S}_{k-1,N} = \sum_{j=1}^{p} X_{k-j,N}^2$. An additional advantage of this weight function, $\kappa(u_0, \mathcal{X}_{k-1,N})$, is that under Assumption 1, $\sup_{k,N} \mathbb{E}(\frac{X_{k-i,N}^4}{\kappa(u_0, \mathcal{X}_{k-1,N})^2}) < \infty$ and $\sup_{k,N} \mathbb{E}(\frac{X_{k-i,N}^{8+2\nu}}{\kappa(u_0, \mathcal{X}_{k-1,N})^{4+\nu}}) < \infty$ are immediately satisfied. Furthermore, $|\kappa(u, \underline{x}) - \kappa(u, \underline{y})| \leq K \sum_{j=1}^{p} |x_j - y_j|$, thus $|x_i/\kappa(u, \underline{x}) - y_i/\kappa(u, \underline{y})| \leq K \sum_{j=1}^{p} |x_j - y_j|$. Therefore, all the conditions in Proposition 3 hold.

3.4. *The two-stage kernel-NLS estimator.* We use $\hat{\mu}_{t_0,N}$ as an estimator of $\mu(u_0)$ (see Lemma A.1 in the Appendix), where

$$\hat{\mu}_{t_0,N} = \sum_{k=1}^{N} \frac{1}{bN} W\left(\frac{t_0 - k}{bN}\right) X_{k,N}^2. \tag{15}$$

We use this to define the two-stage kernel-NLS estimator of the tvARCH($p$) parameters.

*The two-stage scheme*:

(i) Evaluate $\hat{\mu}_{t_0,N}$, given in (15), which is an estimator of $\mu(u_0)$.

(ii) Let $\underline{\hat{a}}_{t_0,N} = \{\tilde{\mathcal{R}}_{t_0,N}\}^{-1} \underline{\tilde{r}}_{t_0,N}$ with $\mathcal{S}_{k-1,N} = \sum_{j=1}^{p} X_{k-j,N}^2$, $\kappa_{t_0,N}(\mathcal{S}_{k-1,N}) = (\hat{\mu}_{t_0,N} + \mathcal{S}_{k-1,N})$ and

$$\tilde{\mathcal{R}}_{t_0,N} = \sum_{k=p+1}^{N} \frac{1}{bN} W\left(\frac{t_0 - k}{bN}\right) \frac{\mathcal{X}_{k-1,N} \mathcal{X}_{k-1,N}^T}{\kappa_{t_0,N}(\mathcal{S}_{k-1,N})^2},$$

$$\tilde{\underline{r}}_{t_0,N} = \sum_{k=p+1}^{N} \frac{1}{bN} W\left(\frac{t_0 - k}{bN}\right) \frac{X_{k,N}^2 \mathcal{X}_{k-1,N}}{\kappa_{t_0,N}(\mathcal{S}_{k-1,N})^2}. \tag{16}$$

If $|u_0 - t_0/N| < 1/N$, we use $\underline{\hat{a}}_{t_0,N}$ as an estimator of $\underline{a}(u_0)$. We call $\underline{\hat{a}}_{t_0,N}$ the two-stage kernel-NLS estimator.

3.5. *Asymptotic properties of the two-stage kernel-NLS estimator.* We derive the asymptotic sampling properties of $\underline{\hat{a}}_{t_0,N}$. [We note that because in the first stage we need to estimate the weight function $\kappa(u_0, \mathcal{X}_{k-1}) = \mu(u_0) + \mathcal{S}_{k-1,N}$, we require the additional mild Assumption 1(iv), which we use to obtain a rate of convergence for $|\hat{\mu}_{t_0,N} - \mu(u_0)|$.]

In the following proposition we obtain consistency and asymptotic normality of $\underline{\tilde{a}}_{t_0,N}$.

PROPOSITION 4. *Suppose* $\{X_{t,N}\}_t$ *is a* tvARCH$(p)$ *process which satisfies Assumption* 1(i), (ii), *and let* $\hat{\mu}_{t_0,N}, \underline{\tilde{a}}_{t_0,N}, \mathcal{A}_t^{(\mu)}(u)$ *and* $\mathcal{D}_t^{(\mu)}(u)$ *be defined as in* (15), *the two stage scheme and* (12), *respectively. Further, let* $\mu(u) = \mathbb{E}(\tilde{X}_t^2(u))$, *and suppose* $|u_0 - t_0/N| < 1/N$.

(i) *Then we have* $\underline{\tilde{a}}_{t_0,N} \xrightarrow{\mathcal{P}} \underline{a}(u_0)$, *with* $b \to 0$, $bN \to \infty$ *as* $N \to \infty$.

(ii) *If in addition we assume that Assumption* 1(iii), (iv) *holds, then we have*

$$
(17) \quad \begin{aligned}
&\sqrt{bN}(\underline{\tilde{a}}_{t_0,N} - \underline{a}(u_0)) + \tfrac{1}{2}\sqrt{bN}\{\mathbb{E}[\mathcal{A}_t^{(\mu)}(u_0)]\}^{-1}\nabla\tilde{\mathcal{B}}_{t_0,N}(\underline{a}(u_0)) \\
&\quad \xrightarrow{\mathcal{D}} \mathcal{N}(0, w_2\mu_4\{\mathbb{E}[\mathcal{A}_t^{(\mu)}(u_0)]\}^{-1}\mathbb{E}[\mathcal{D}_t^{(\mu)}(u_0)]\{\mathbb{E}[\mathcal{A}_t^{(\mu)}(u_0)]\}^{-1}),
\end{aligned}
$$

*where* $\nabla\tilde{\mathcal{B}}_{t_0,N}(\underline{a}(u_0)) = O_p(b^\beta)$ *and* $w_2$ *and* $\mu_4$ *are defined as in Proposition* 3, *with* $b \to 0$, $bN \to \infty$ *as* $N \to \infty$.

Comparing the two-stage kernel-NLS estimator with the kernel-QML estimator in Dahlhaus and Subba Rao (2006), it is easily seen that they both have the same rate of convergence.

REMARK 1 (An asymptotically optimal estimator). We recall that the oracle estimator asymptotically has the same variance as the kernel-QML estimator, but in practice the oracle weight is never known. However, the two-stage kernel-NLS estimator can be used as the basis of an estimate of the oracle weight. In other words, using the two-stage kernel-NLS estimator, we define the weight function $\hat{\sigma}_{k,N}^2(u_0) = \tilde{a}_{t_0,N}(0) + \sum_{j=1}^p \tilde{a}_{t_0,N}(j)X_{k-j,N}^2$, where $\underline{\tilde{a}}_{t_0,N} = (\tilde{a}_{t_0,N}(0), \ldots, \tilde{a}_{t_0,N}(p))$. Then, we use $\underline{\breve{a}}_{t_0,N}$ as an estimator of $\underline{a}(u_0)$, where $\underline{\breve{a}}_{t_0,N} = \{\breve{\mathcal{R}}_{t_0,N}\}^{-1}\breve{\underline{r}}_{t_0,N}$, and $\breve{\mathcal{R}}_{t_0,N}$ and $\breve{\underline{r}}_{t_0,N}$ are defined in the same way as $\tilde{\mathcal{R}}_{t_0,N}$ and $\tilde{\underline{r}}_{t_0,N}$, with $\hat{\sigma}_{t,N}^2(u_0)$ replacing $(\hat{\mu}_{t_0,N} + \sum_{j=1}^p X_{t-j,N}^2)$. The asymptotic sampling results can be derived using a similar proof to Proposition 4. More precisely, if Assumption 1 holds, $b^\beta\sqrt{bN} \to 0$, and $a_j(u_0) > 0$ for all $j$, then we have

$$
(18) \quad \sqrt{bN}(\underline{\breve{a}}_{t_0,N} - \underline{a}(u_0)) \xrightarrow{\mathcal{D}} \mathcal{N}(0, w_2\mu_4\{\mathbb{E}[\Sigma_t(u_0)]\}^{-1}).
$$

In other words, by using the two-stage kernel-NLS estimator, we are able to estimate the oracle weight sufficiently well for the parameter to have the same asymptotic variance as the kernel-QML estimator. We note that, similarly to the kernel-QML estimator, we require that $\inf_j a_j(u) > 0$. However, it is suggested in the examples in Section 6 that over large periods of time some of the higher-order parameters should be zero. This renders

*Ratios of Mean Absolute Errors of two-stage NLS and QML estimators, averaged over 100 simulated sample paths, for stationary* ARCH(2) *estimation with Gaussian errors* $Z_t$ *and* $(a_0, a_1, a_2) = (1, 0.6, 0.3)$. *Sample sizes vary from* $N = 15$ *to* $N = 250$

|       | $N = 15$ | $N = 30$ | $N = 60$ | $N = 100$ | $N = 150$ | $N = 250$ |
|-------|----------|----------|----------|-----------|-----------|-----------|
| $a_0$ | 0.59     | 0.69     | 0.91     | 1.04      | 0.96      | 1.28      |
| $a_1$ | 0.84     | 0.73     | 0.97     | 0.97      | 1.10      | 1.11      |
| $a_2$ | 0.64     | 0.68     | 0.86     | 0.98      | 0.94      | 1.08      |

the assumption $\inf_j a_j(u) > 0$ rather unrealistic. Furthermore, to estimate $\breve{\underline{a}}_{t_0,N}$, we require an additional stage of computation, which significantly increases computation time in tasks such as cross-validatory bandwidth choice or evaluation of bootstrap confidence intervals. Also, small sample evidence suggests that the performance of the estimators $\underline{\tilde{a}}_{t_0,N}$ and $\breve{\underline{a}}_{t_0,N}$ is similar. For this reason, in the rest of this paper, we focus on $\underline{\tilde{a}}_{t_0,N}$, though our results can be generalized to $\breve{\underline{a}}_{t_0,N}$.

3.6. *Comparison of two-stage kernel-NLS and kernel-QML estimators for small samples.* As mentioned earlier, in a nonstationary setting, it is essential for any estimator of tvARCH($p$) parameters to perform well for small sample sizes. We now briefly describe the outcome of an extensive simulation study aimed at comparing the performance of the two-stage NLS and QML estimators on short stretches of stationary ARCH(2) data. We have tested the two estimators for Gaussian, Laplace and Student-$t$ errors $Z_t$, and for various points of the parameter space $(a_0, a_1, a_2)$. The two-stage NLS estimator significantly outperformed the QML estimator for very small sample sizes in almost all of the cases. More complicated patterns emerged for sample sizes of about 150 and larger, where the performance depended on the particular point of the parameter space. However, the two-stage NLS estimator was never found to perform much worse than the QML estimator. We also found the two-stage NLS estimator to be significantly faster than the QML estimator as it did not involve an iterative optimization procedure.

As an example, Table 1 shows the ratios of the mean absolute errors of the two-stage NLS and QML estimators, averaged over 100 simulated sample paths, for the following parameter configuration: $(a_0, a_1, a_2) = (1, 0.6, 0.3)$. The errors $Z_t$ are Gaussian. The above point of the parameter space is "typical" in the sense that it lies in the interior of the parameter space (and thus is suitable for QML estimation which requires $a_1, a_2 > 0$) and that $a_1 > a_2$ as expected in a real-data setting. Also, it is interesting in that $a_1 + a_2 > 1/\sqrt{3}$ and thus the classical (nonnormalized) least-squares estimator, corresponding to $\kappa(\cdot) \equiv 1$, would not be consistent in this setup.

**4. A cross-validation method for bandwidth selection and implementation.** In this section, we propose a data-driven method for selecting the bandwidth of the two-stage kernel-NLS estimator.

4.1. *The cross-validation bandwidth estimator.* Several cross-validation methods in nonparametric statistics consider the distance between an observation and a predictor of that observation given neighboring observations. For example, Hart ([1996](#)) used a cross-validation method based on the best linear predictor of $Y_t$ given the past to select the bandwidth of a kernel smoother, where $Y_t$ was a nonparametric function plus correlated noise. The methodology we propose is based on the best linear predictor of $X_{t,N}^2$ given the past, which is $a_0(\frac{t}{N}) + \sum_{j=1}^p a_j(\frac{t}{N})X_{t-j,N}^2$.

We estimate the parameters $\{a_j(t/N)\}_j$ using the localized two-stage kernel-NLS method but omit the observation $X_{t,N}^2$ in the estimation. More precisely, we use $\underline{\hat{a}}_{t,N}^{-t}(b) = (\hat{a}_0^{-t}(b), \ldots, \hat{a}_p^{-t}(b))$ as an estimator of $\{a_j(t/N)\}_j$, where

$$(19) \qquad \underline{\hat{a}}_{t,N}^{-t}(b) = \{\mathcal{R}_{t,N}^{-t}(b)\}^{-1}\underline{r}_{t,N}^{-t}(b),$$

with

$$\mathcal{R}_{t,N}^{-t}(b) = \sum_{\substack{k=p+1 \\ k \neq t, \ldots, t+p}}^N \frac{1}{bN}W\left(\frac{t-k}{bN}\right)\frac{\mathcal{X}_{k-1,N}\mathcal{X}_{k-1,N}^T}{(\hat{\mu}_{t,N} + \mathcal{S}_{k-1,N})^2},$$

$$\underline{r}_{t,N}^{-t}(b) = \sum_{\substack{k=p+1 \\ k \neq t, \ldots, t+p}}^N \frac{1}{bN}W\left(\frac{t-k}{bN}\right)\frac{X_{k,N}^2\mathcal{X}_{k-1,N}}{(\hat{\mu}_{t,N} + \mathcal{S}_{k-1,N})^2}.$$

By using $\underline{\hat{a}}_{t,N}^{-t}(b)$, the squared error in predicting $X_{t,N}^2$ is given by $(X_{t,N}^2 - \hat{a}_0^{-t}(b) - \sum_{j=1}^p \hat{a}_j^{-t}(b)X_{t-j,N}^2)^2$.

To reduce the complexity, we suggest only evaluating the cross-validation criterion on a subsample of the observations. Let $h$ be such that $h \to \infty$, $N/h \to \infty$ as $N \to \infty$ (in practice $h \gg p$). We implement the cross-validation criterion on only the subsampled observations $\{X_{kh,N} : k = 1, \ldots, N/h\}$. In other words, let $\underline{\hat{a}}_{kh,N}^{-kh}(b) = (\hat{a}_0^{-kh}(b), \ldots, \hat{a}_p^{-kh}(b))$ be the estimator defined in [(19)](#) and by normalizing the squared error with the term $(\hat{\mu}_{kh,N} + \sum_{j=1}^p X_{kh-j,N}^2)^2$, we define the following cross-validation criterion

$$(20) \qquad \mathcal{G}_{N,h}(b) = \frac{h}{N}\sum_{k=1}^{N/h}\frac{(X_{kh,N}^2 - \hat{a}_0^{-kh}(b) - \sum_{j=1}^p \hat{a}_j^{-kh}(b)X_{kh-j,N}^2)^2}{(\hat{\mu}_{kh,N} + \sum_{j=1}^p X_{kh-j,N}^2)^2}.$$

We then use $\hat{b}_{\text{opt}}^h$ as the optimal bandwidth, where $\hat{b}_{\text{opt}}^h = \arg\min_b \mathcal{G}_{N,h}(b)$. Using similar arguments to those in Hart ([1996](#)), asymptotically, one can show that $\mathcal{G}_{N,h}(b)$ is equivalent to the mean-squared error $\tilde{\mathcal{G}}_{N,h}(b)$, where

$$(21) \quad \tilde{\mathcal{G}}_{N,h}(b) = \frac{h}{N} \sum_{k=1}^{N/h} \mathbb{E}\left\{ \frac{(X_{kh,N}^2 - \hat{a}_0^{-kh}(b) - \sum_{j=1}^p \hat{a}_j^{-kh}(b) X_{kh-j,N}^2)^2}{(\hat{\mu}_{kh,N} + \sum_{j=1}^p X_{kh-j,N}^2)^2} \right\}.$$

It follows that $\hat{b}_{\mathrm{opt}}^h$ is an estimator of $b_{\mathrm{opt}}$, where $b_{\mathrm{opt}} = \arg\min_b \tilde{\mathcal{G}}_{N,h}(b)$. $\tilde{\mathcal{G}}_{N,h}(b)$ is minimized if $\hat{\underline{a}}_{kh,N}^{-t}(b) = \underline{a}(kh/N)$ and in that case it is asymptotically equal to

$$(22) \quad \int_0^1 \mathbb{E}\left\{ \frac{(Z_0^2 - 1)^2 \sigma_0^2(u)}{[\mu(u) + \sum_{j=1}^p X_{-j}^2(u)]^2} \right\} du.$$

Therefore, $\hat{b}_{\mathrm{opt}}^h$ is such that $\hat{\underline{a}}_{t,N}^{-t}(\hat{b}_{\mathrm{opt}}^h)$ is close to $\underline{a}(t/N)$.

It is straightforward to show that the computational complexity of this algorithm is $O(B\frac{N}{h}N\log N)$, where $B$ is the cardinality of the set of bandwidths tested for the minimum of the cross-validation criterion. We note that the above rate is unattainable for the kernel-QML estimator due to its iterative character.

4.2. *An illustrative example.* We illustrate the performance of the proposed cross-validation criterion by an interesting example of a tvARCH(1) process for which the parameters $a_0(\cdot)$ and $a_1(\cdot)$ vary over time but the asymptotic unconditional variance $\mathbb{E}(\tilde{X}_t^2(u)) = a_0(u)/(1 - a_1(u))$ remains constant. This means that sample paths of $\{X_{t,N}\}_t$ will invariably appear stationary on visual inspection, and that more sophisticated techniques are needed to detect the nonstationarity.

The left-hand plot in Figure 1 shows a sample path of length 1024, simulated from the above process using standard Gaussian errors. The true time-varying parameters $a_0(\cdot)$ and $a_1(\cdot)$ are displayed as dotted lines in the middle and right-hand plots, respectively. In the estimation procedure, we used the Parzen kernel (a convolution of the rectangular and triangular kernels) and, for simplicity, set $\hat{\mu}_{t,N}$ to be the sample mean of $\{X_{t,N}^2\}_t$. To estimate a suitable bandwidth, we applied the proposed cross-validation procedure described above with $h = 10$ (empirically, we have found that for data of length of order 1000, the value $h = 10$ offers a good compromise between speed and accuracy of our method). We examined the value of the cross-validation criterion over a regular grid of bandwidths between 0 and 1, and obtained the optimal bandwidth as $\hat{b}_{\mathrm{opt}}^h = 0.132$.

The resulting parameter estimates are shown in the middle and right-hand plots of Figure 1 as solid lines. While we can clearly observe a degree of bias due to the small sample sizes involved in the estimation, it is reassuring to see that the resulting estimates correctly trace the shape of the underlying parameters. Denoting the empirical residuals from the fit by $\hat{Z}_t$, the $p$-value of the Kolmogorov–Smirnov test for Gaussianity of $\hat{Z}_t$ was 0.08, and the $p$-values of the Ljung–Box test for lack of serial correlation in $\hat{Z}_t$, $|\hat{Z}_t|$ and $\hat{Z}_t^2$ were 0.71, 0.33 and 0.58, respectively.

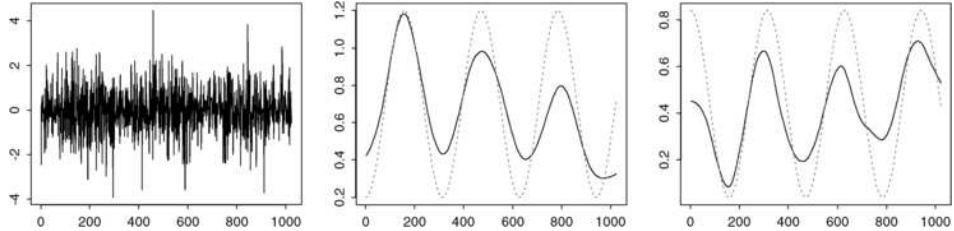

Fig. 1. *Dotted lines in the middle and right plots: the true time-varying parameters $a_0(u)$ and $a_1(u)$, respectively. The left plot: a sample path from the model, with Gaussian errors. Solid lines in the middle and right plots: the corresponding estimates. See Section 4.2 for details.*

## 5. Constructing bootstrap pointwise confidence intervals.

In parameter estimation of linear time series, bootstrap methods are often used to obtain a good finite sample approximation of the distribution of the parameter estimators. Schemes based on estimating the residuals are often used [Franke and Kreiss (1992)]. Inspired by these methods, we propose a bootstrap scheme for the tvARCH($p$) process, which we use to construct pointwise confidence intervals for the two-stage kernel-NLS estimator. The main idea of the scheme is to use the two-stage kernel-NLS estimator to estimate the residuals. We construct the empirical distributions from the estimated residuals, sample from it and use this to construct the bootstrap tvARCH($p$) sample. We show that the distribution of the two-stage kernel-NLS estimator using the bootstrap tvARCH($p$) sample and the "true" tvARCH($p$) estimator asymptotically coincide. We mention that the scheme and the asymptotic results derived here are also of independent interest and can be used to bootstrap stationary ARCH($p$) processes [for a recent review on resampling and subsampling financial time series in the stationary context, see Paparoditis and Politis (2007)]. We emphasize that unlike the kernel-QML estimator, this computer-intensive procedure is feasible for the kernel-NLS estimator due to its rapid computability.

Let $\underline{\tilde{a}}_{t_0,N} = (\tilde{a}_{t_0,N}(0), \ldots, \tilde{a}_{t_0,N}(p))$. We first note that Assumption 1(i) is usually imposed in the tvARCH framework because it guarantees that almost surely every realization of the resulting process is bounded. When the sum of the coefficients is greater than one, the corresponding process is unstable. The following residual bootstrap scheme constructs the tvARCH($p$) process from estimates of the residuals and the parameter estimators. Despite $\underline{\tilde{a}}_{t_0,N} \xrightarrow{\mathcal{P}} \underline{a}(u_0)$, it is not necessarily true that the sum of the parameter estimates satisfies $\sum_{j=1}^{p} \tilde{a}_{t_0,N}(j) < 1$. To overcome this, we now define a very slight modification of the two-stage kernel-NLS estimator which guarantees that this sum is less than one. Let $\underline{\bar{a}}_{t_0,N} = (\bar{a}_{t_0,N}(0), \ldots, \bar{a}_{t_0,N}(p))$, where $\bar{a}_{t_0,N}(0) = \tilde{a}_{t_0,N}(0)$ and, for $j > 1$,

$$(23) \quad \bar{a}_{t_0,N}(j) = \begin{cases} \tilde{a}_{t_0,N}(j), & \text{if } \sum_{j=1}^{p} \tilde{a}_{t_0,N}(j) \leq 1 - \delta, \\[2mm] (1-\delta)\dfrac{\tilde{a}_{t_0,N}(j)}{\sum_{j=1}^{p} \tilde{a}_{t_0,N}(j)}, & \text{if } \sum_{j=1}^{p} \tilde{a}_{t_0,N}(j) > 1 - \delta. \end{cases}$$

Since $\underline{\tilde{a}}_{t_0,N} \xrightarrow{\mathcal{P}} \underline{a}(u_0)$ and $\sum_{j=1}^{p} a_j(u) \leq 1 - \delta$ [Assumption 1(i)], it is straightforward to see that $\underline{\bar{a}}_{t_0,N} \xrightarrow{\mathcal{P}} \underline{a}(u_0)$ and $\sum_{j=1}^{p} \bar{a}_{t_0,N}(j) \leq 1 - \delta$.

*The residual bootstrap of the* tvARCH($p$) *process:*

(i) If $k \in [t_0 - bN, t_0 + bN - 1]$, using the parameter estimators construct residuals

$$\tilde{Z}_k^2 = \frac{X_{k,N}^2}{\bar{a}_{t_0,N}(0) + \sum_{j=1}^{p} \tilde{a}_{t_0,N}(j) X_{k-j,N}^2}.$$

(ii) Define $\hat{Z}_t^2 = \tilde{Z}_t^2 - \frac{1}{2bN} \sum_{k=t_0-bN}^{t_0+bN-1} \tilde{Z}_k^2 + 1$ and consider the empirical distribution function

$$\hat{F}_{t_0,N}(x) = \frac{1}{bN} \sum_{k=t_0-bN}^{t_0+bN-1} \mathbb{I}_{(-\infty,x]}(\hat{Z}_k^2),$$

where $\mathbb{I}_A(y) = 1$ if $y \in A$, 0 otherwise. It is worth mentioning that we use $\hat{Z}_t^2$ rather than $\tilde{Z}_t^2$ since we have $\mathbb{E}(\hat{Z}_t^2) = \int z \hat{F}_{t_0,N}(dz) = 1$. (This result is used in Proposition 6 in the Appendix.)

Set $X_t^{+2}(u_0) = 0$ for $t \leq 0$. For $1 \leq t \leq t_0 + bN/2$, sample from the distribution function $\hat{F}_{t_0,N}(x)$, to obtain the sample $\{Z_t^{+2}\}_t$. Use this to construct the bootstrap sample

$$X_t^{+2}(u_0) = \sigma_t^{+2}(u_0) Z_t^{+2}, \qquad \sigma_t^{+2}(u_0) = \bar{a}_{t_0,N}(0) + \sum_{j=1}^{p} \bar{a}_{t_0,N}(j) X_{t-j}^{+2}(u_0).$$

We note that by estimating the residuals from $[t_0 - bN, t_0 + bN - 1]$, the distribution of $X_t^{+2}(u_0)$ will be suitably close to the stationary approximation $X_t(u_0)$ when $t \in [t_0 - bN/2, t_0 + bN/2 - 1]$, this allows us to obtain the sampling properties of the bootstrap estimator.

(iii) Define the bootstrap estimator

$$(24) \qquad \underline{\hat{a}}_{t_0,N}^{+} = \{\mathcal{R}_{t_0,N}^{+}\}^{-1} \underline{r}_{t_0,N}^{+},$$

where $\mathcal{X}_{t-1}(u_0)^+ = (1, X_{t-1}^{+2}(u_0), \ldots, X_{t-p}^{+2}(u_0))^T$ and

$$\mathcal{R}_{t_0,N}^+ = \sum_{k=p+1}^{N} \frac{1}{bN} W\left(\frac{t_0-k}{bN}\right) \frac{\mathcal{X}_{k-1}(u_0)^+ \mathcal{X}_{k-1}(u_0)^{+T}}{(\hat{\mu}_{t_0,N} + \sum_{j=1}^{k} X_{k-j}^{+2}(u_0))^2},$$

$$\mathcal{L}_{t_0,N}^+ = \sum_{k=p+1}^{N} \frac{1}{bN} W\left(\frac{t_0-k}{bN}\right) \frac{X_k^{+2}(u_0) \mathcal{X}_{k-1}(u_0)^{+T}}{(\hat{\mu}_{t_0,N} + \sum_{j=1}^{k} X_{k-j}^{+2}(u_0))^2}.$$

We observe that in steps (i), (ii) of the bootstrap scheme we are constructing the bootstrap sample $\{X_t^{+2}(u_0)\}_t$ whose distribution should emulate the distribution of the stationary approximation $\{\tilde{X}_t^2(u_0)\}_t$. In step (iii) of the bootstrap scheme we are constructing the bootstrap estimator $\hat{\underline{a}}_{t_0,N}^+$ from the bootstrap samples. We note that we have bootstrapped the stationary approximation $\tilde{X}_t^2(u_0)$ since the limiting distribution of $\hat{\underline{a}}_{t_0,N}$ is derived using the stationary approximation.

We now show that the distributions of $\sqrt{bN}\{\hat{\underline{a}}_{t_0,N}^+ - \bar{\underline{a}}_{t_0,N}\}$ and $\sqrt{bN}\{\hat{\underline{a}}_{t_0,N} - \underline{a}(u_0)\}$ asymptotically coincide.

PROPOSITION 5. *Suppose Assumption* 1 *holds, and suppose either* $\inf_j a_j(u_0) > 0$ *or* $\mathbb{E}(Z_t^4)^{1/2} \sup_u [\sum_{j=1}^{p} a_j(u)] < 1 - \delta$ *[which implies* $\sup_k \mathbb{E}(X_{k,N}^4) < \infty$*]. Let* $\bar{\underline{a}}_{t_0,N}$ *and* $\hat{\underline{a}}_{t_0,N}^+$ *be defined as in* (23) *and* (24), *respectively, and let* $b^\beta \sqrt{bN} \to 0$. *If* $|u_0 - t_0/N| < 1/N$, *then we have*

$$\sqrt{bN}(\hat{\underline{a}}_{t_0,N}^+ - \bar{\underline{a}}_{t_0,N})$$
$$\xrightarrow{\mathcal{D}} \mathcal{N}(0, w_2 \mu_4 \{\mathbb{E}[\mathcal{A}_t^{(\mu)}(u_0)]\}^{-1} \{\mathbb{E}[\mathcal{D}_t^{(\mu)}(u_0)]\} \{\mathbb{E}[\mathcal{A}_t^{(\mu)}(u_0)]\}^{-1}),$$

*with* $b \to 0$, $bN \to \infty$ *as* $N \to \infty$.

Comparing the results in Propositions 4(ii) and Propositions 5 we see if $b^\beta \sqrt{bN} \to 0$, then, asymptotically, the distributions of $(\hat{\underline{a}}_{t_0,N}^+ - \bar{\underline{a}}_{t_0,N})$ and $(\hat{\underline{a}}_{t_0,N} - \underline{a}(u_0))$ are the same.

**6. Volatility estimation: real data examples.** The datasets analyzed in this and the following section fall into two categories:

1. Logged and differenced daily exchange rates between USD and a number of other currencies running from January 1, 1990 to December 31, 1999: the data are available from the US Federal Reserve website: **www.federalreserve.gov/releases/h10/Hist/default1999.htm**. We use the following acronyms: CHF (Switzerland Franc), GBP (United Kingdom Pound), HKD (Hong Kong Dollar), JPY (Japan Yen), NOK (Norway Kroner), NZD (New Zealand Dollar), SEK (Sweden Kronor), TWD (Taiwan New Dollar).

2. Logged and differenced daily closing values of the NIKKEI, FTSE, S and P500 and DAX indices, measured between a date in 1996 (exact dates vary) and April 29, 2005: the data are available from: www.bossa.pl/notowania/daneatech/metastock/.

The lengths $N$ of each dataset vary but oscillate around 2500. In this section, we exhibit the estimation performance of the two-stage kernel-NLS estimator on the USD/GBP exchange rate and FTSE series. We examine the cases $p = 0, 1, 2$ and use the Parzen kernel with bandwidths selected by the cross-validation algorithm of Section 4.2.

The left column in Figure 2 shows the results for USD/GBP. The top plot shows the data, the next one down shows the estimates of $a_0(\cdot)$ for $p = 0$ (dashed line), $p = 1$ (dotted line) and $p = 2$ (solid line), the one below displays the positive parts of the estimates of $a_1(\cdot)$ for $p = 1$ (dotted) and $p = 2$ (solid), and the bottom plot shows the positive part of the estimate of $a_2(\cdot)$ for $p = 2$. Note that the negative values arise since our estimator is not guaranteed to be nonnegative. The right column shows the corresponding quantities for the FTSE data. It is interesting to observe that in both cases, the shapes of the estimated time-varying parameters are similar for different values of $p$.

The goodness of fit for each choice of $p = 0, 1, 2$ is assessed in Table 2. In each case, $\hat{Z}_t$ denotes the sequence of empirical residuals from the given fit. For the USD/GBP data, the best fit is obtained for $p = 1$. For the FTSE data, it is less clear which order gives the best fit but the Ljung–Box (L–B) $p$-value for $|\hat{Z}_t|$ is the highest for $p = 0$ and thus it seems to be the preferred option, which is further confirmed by the visual inspection of the sample autocorrelation function of $|\hat{Z}_t|$ in the three cases. In both cases, the empirical residuals are negatively skewed, and in the case of USD/GBP they are also heavy-tailed.

We conclude this section by constructing bootstrap pointwise confidence intervals for the estimated parameters, using the algorithm detailed in Section 5. Note that our central limit theorem (CLT) of Proposition 4 could be used for the same purpose, but this would require pre-estimation of a number of quantities, which we wanted to avoid. We base our bootstrap pointwise confidence intervals on 100 bootstrap samples. For clarity, we only display confidence intervals for the "preferred" orders $p$: that is, for $p = 1$ in the case of the USD/GBP data, and $p = 0$ in the case of the FTSE series. These are shown in Figure 3.

It is interesting to note that the pointwise confidence intervals for the "nonlinearity" parameter $a_1(\cdot)$ in the USD/GBP series are relatively wide and that the parameter can be viewed as only insignificantly different from zero most (but not all) of the time. On the other hand, there exist time intervals where the parameter significantly deviates from zero. This further confirms the observation made earlier that the order $p = 0$ is an inferior modeling choice for this series and that the order $p = 1$ is preferred.

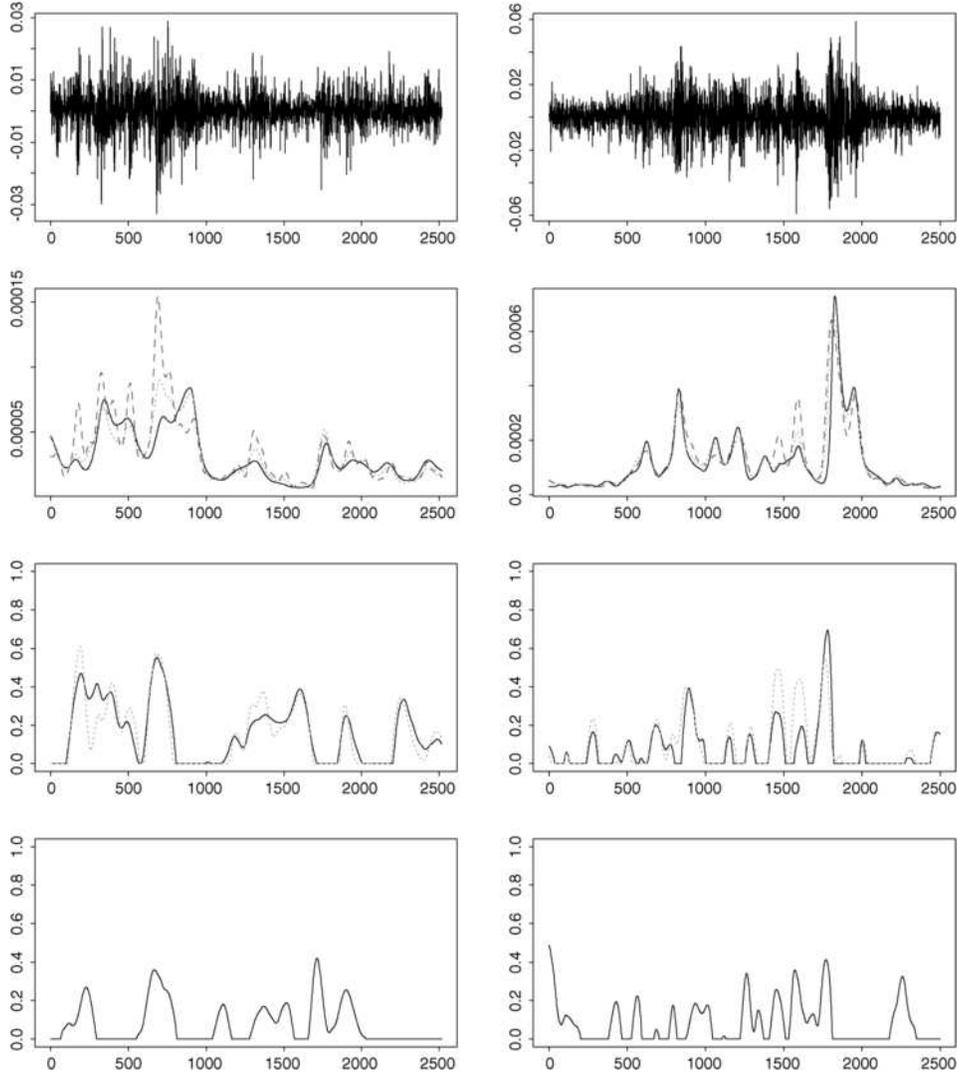

Fig. 2. *Left (right) column: USD/GBP (FTSE) series and the corresponding estimation results. See Section 6 for details.*

**7. Volatility forecasting: real data examples.** In this section, we describe a numerical study whereby the long-term volatility forecasting ability of the tvARCH($p$) process is compared to that of the stationary GARCH($1, 1$) and EGARCH($1, 1$) processes with standard Gaussian errors. We compute the forecasts of the tvARCH($p$) process as follows: we use the available data to

TABLE 2

*The values of bandwidth selected by cross-validation, the p-values of the L–B test for white noise for $\hat{Z}_t$, $|\hat{Z}_t|$, $\hat{Z}_t^2$, and the sample skewness and kurtosis coefficients for $\hat{Z}_t$ for the USD/GBP and FTSE data sets. The boxed value means p-value is below 0.05*

|  | USD/GBP | | | FTSE | | |
|---|---|---|---|---|---|---|
|  | $p = 0$ | $p = 1$ | $p = 2$ | $p = 0$ | $p = 1$ | $p = 2$ |
| Bandwidth | 0.02 | 0.032 | 0.04 | 0.024 | 0.028 | 0.028 |
| L–B $P$-value for $\hat{Z}_t$ | 0.83 | 0.83 | 0.82 | 0.15 | 0.20 | 0.30 |
| L–B $P$-value for $|\hat{Z}_t|$ | 0.17 | 0.71 | 0.03 | 0.10 | 0.07 | 0.07 |
| L–B $P$-value for $\hat{Z}_t^2$ | 0.09 | 0.79 | 0.26 | 0.13 | 0.35 | 0.52 |
| Skewness of $\hat{Z}_t$ | −0.05 | −0.09 | −0.08 | −0.13 | −0.15 | −0.16 |
| Kurtosis of $\hat{Z}_t$ | 0.7 | 0.92 | 1.24 | −0.01 | 0.06 | 0.15 |

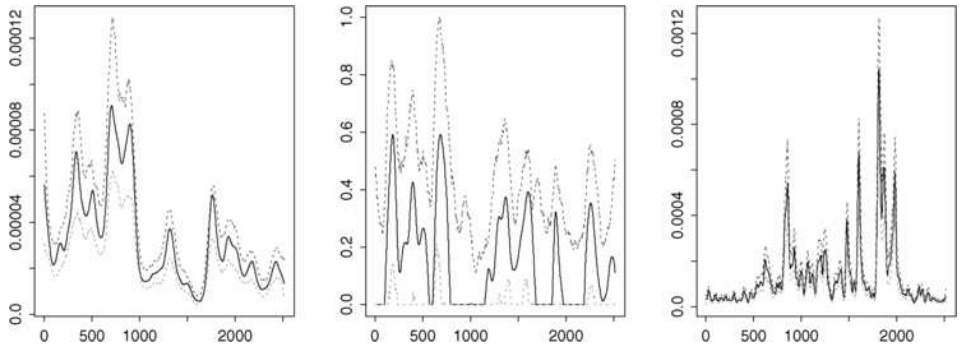

FIG. 3. *Solid lines from left to right: estimates of $a_0(\cdot)$ for USD/GBP, $a_1(\cdot)$ for USD/GBP, and $a_0(\cdot)$ for FTSE. Dashed lines: the corresponding 80% symmetric bootstrap pointwise confidence intervals.*

estimate the tvARCH($p$) parameters, and then forecast into the future using the "last" estimated parameter values, that is, those corresponding to the right edge of the observed data. For a rectangular kernel with span $m$, this strategy leads to the following algorithm: (a) treat the last $m$ data points as if they came from a stationary ARCH($p$) process, (b) estimate the stationary ARCH($p$) parameters on this segment (via the two-stage NLS scheme), and (c) forecast into the future as in the classical stationary ARCH($p$) forecasting theory [for the latter, see, e.g., Bera and Higgins (1993)].

We denote the mean-square-optimal $h$-step-ahead volatility forecasts at time $t$, obtained via the above algorithm, by $\sigma_{t|t+h}^{2,\text{tvARCH}(p)}$. Note that to obtain the analogous quantities, $\sigma_{t|t+h}^{2,\text{GARCH}(1,1)}$ and $\sigma_{t|t+h}^{2,\text{EGARCH}(1,1)}$, for the stationary GARCH(1, 1) and EGARCH(1, 1) processes, we always use the entire available dataset, and not only the last $m$ observations.

To test the forecasting ability of the various models, we use the exchange rate and stock index datasets listed in Section 6. For the tvARCH($p$) process, we take $p = 0, 1, 2$, and use the forecasting procedure described above with a rectangular kernel, over a grid of span values $m = 50, 100, \ldots, 500$. Note that the tvARCH(0) process has the simple form $X_{t,N} = a_0^{1/2}(t/N)Z_t$ and is also considered by Stărică and Granger (2005). We select the span by a "forward validation" procedure, that is, choose the value of $m$ that yields the minimum out-of-sample prediction error AMSE defined below.

For the stationary (E)GARCH(1,1) prediction, we use the standard S-Plus `garch` and `predict` routines. The stationary (E)GARCH(1,1) parameters are re-estimated for each $t$.

For each $t = 1000, \ldots, N - 250$, we compute the quantities

$$\overline{\sigma}^{2,\text{model}}_{t|t+250} = \sum_{h=1}^{250} \sigma^{2,\text{model}}_{t|t+h},$$

where "model" is one of: tvARCH(0), tvARCH(1), tvARCH(2), GARCH(1,1), and EGARCH(1,1), and compare them to the "realized" volatility

$$\overline{X}^2_{t|t+250} = \sum_{h=1}^{250} X^2_{t+h},$$

using the scaled aggregated mean square error (AMSE)

$$R^{\text{model}}_{250,1000,N} = \sum_{t=1000}^{N-250} (\overline{\sigma}^{2,\text{model}}_{t|t+250} - \overline{X}^2_{t|t+250})^2,$$

where the scaling is by the factor of $1/(N - 1000)$. For a justification of this simulation setup, see Stărică (2003).

Table 3 lists the AMSEs attained by tvARCH(0), tvARCH(1), tvARCH(2), stationary GARCH(1,1) and stationary EGARCH(1, 1) processes: the best results are boxed. The values in brackets indicate the selected span values. The bullets for the USD/TWD and USD/HKD series indicate that the numerical optimizers performing the QML estimation in stationary (E)GARCH(1,1) processes failed to converge at several points of the series and, therefore, we were unable to obtain accurate forecasts. We list below some interesting conclusions from this study.

- In most cases, the selected span values $m$ are similar across orders $p$. These values can be taken as an indication of how "variable" the time-varying parameters are. Exceptions to this rule occur mostly in data sets which are difficult to model, such as the HKD series, which is extremely spiky. For the latter series, more thought is needed on how to model it accurately in the tvARCH($p$) (or indeed any other) framework.

Table 3

*AMSE for long-term forecasts using* tvARCH(0), tvARCH(1), tvARCH(2), *stationary* GARCH(1,1) *and stationary EGARCH(1, 1) processes.* $R^{(\mathrm{E})\mathrm{GARCH}(1,1)}_{250,1000,N}$ *is the better result out of:* $R^{\mathrm{GARCH}(1,1)}_{250,1000,N}$ *and* $R^{\mathrm{EGARCH}(1,1)}_{250,1000,N}$

| Series | Scaling | $R^{(\mathrm{E})\mathrm{GARCH}(1,1)}_{250,1000,N}$ | $R^{\mathrm{tvARCH}(0)}_{250,1000,N}$ | $R^{\mathrm{tvARCH}(1)}_{250,1000,N}$ | $R^{\mathrm{tvARCH}(2)}_{250,1000,N}$ |
|---|---|---|---|---|---|
| CHF | $10^8$ | 2395 | 2371 (500) | $\boxed{2254}$ (500) | 3030 (500) |
| GBP | $10^9$ | 20282 | $\boxed{7660}$ (250) | 9567 (300) | 9230 (300) |
| HKD | $10^{12}$ | • | 230 (150) | 170 (500) | $\boxed{150}$ (100) |
| JPY | $10^8$ | $\boxed{8687}$ | 9713 (350) | 9173 (300) | 9450 (300) |
| NOK | $10^8$ | 1767 | $\boxed{1552}$ (500) | 1875 (250) | 2221 (500) |
| NZD | $10^8$ | 11890 | 5270 (50) | 4976 (100) | $\boxed{4955}$ (150) |
| SEK | $10^9$ | 37720 | $\boxed{6639}$ (250) | 6805 (250) | 7321 (250) |
| TWD | $10^8$ | • | $\boxed{2323}$ (500) | 2372 (500) | 2400 (500) |
| S & P500 | $10^5$ | $\boxed{33}$ | 43 (500) | 43 (500) | 40 (500) |
| FTSE | $10^6$ | $\boxed{516}$ | 860 (500) | 958 (500) | 983 (500) |
| DAX | $10^6$ | 2602 | 4492 (150) | 4483 (500) | 4864 (150) |
| NIKKEI | $10^7$ | $\boxed{2364}$ | 3418 (100) | 3252 (250) | 3432 (250) |

- For the NZD series, it can clearly be seen how "adding more nonlinearity takes away nonstationarity": as $p$ increases, a larger and larger span $m$ is selected, which means that more and more variability in the volatility of the data can be attributed to the nonlinearity, rather than the nonstationarity.

- While the tvARCH($p$) framework seems superior to stationary (E)GARCH(1,1) methodology for the currency exchange data, the opposite is true for the stock indices. This might be indicative of the fact that stock indices are "less nonstationary" than currency exchange series.

We conclude with a heuristic investigation of the quality of our volatility forecasts. Conditioning on the information available up to time $t$, the quantity $\overline{\sigma}^{2,\mathrm{model}}_{t|t+250}$ predicts the variance of the variable $X_t^{(250)} := \sum_{h=1}^{250} X_{t+h}$. By CLT-type arguments, $X_t^{(250)}$ is approximately Gaussian, and thus we assess the quality of the predicted volatility by measuring how often the process $Y_t := X_t^{(250)}/\{\overline{\sigma}^{2,\mathrm{model}}_{t|t+250}\}^{1/2}$ falls into desired confidence intervals for standard Gaussian variables.

However, this is less informative of the quality of the forecasting procedure than one might hope, the reason being that the process $Y_t$ is strongly dependent, so it is not reasonable to expect it to take values outside $(1-\alpha)100\%$ confidence intervals exactly, or approximately, $100\alpha\%$ of the time. Figure 4 shows processes $Y_t$ constructed for the GBP, NZD and SEK series, with the "optimal" forecasting parameters from Table 3 (i.e., those for which the results are boxed). For $\alpha = 0.05$, the coverages are, respectively, 100%, 79% and 95%. If the dependence in $Y_t$ were weaker, we would expect the three coverages to be closer to 95%, provided the forecasting procedure was "adequate." However, here, the strong dependence in $Y_t$ causes the variance of the coverage percentages to be high.

Nonetheless, it is reassuring to note that on average, across the datasets, we do obtain the correct coverage of around 95%. To see this, let us consider the series for which our forecasting procedure is satisfactory [i.e., those for which it outperforms (E)GARCH(1,1) processes], bar the two series: HKD and TWD, which are extremely spiky and thus difficult to model and forecast. These are: CHF, GBP, NOK, NZD, SEK. Table 4 shows the coverages for the five series. The average coverage is 94.2%, which is very close to the ideal coverage of 95%. Averaging across all series, excluding HKD and TWD, we obtain a coverage of 95.7%.

## APPENDIX: AUXILIARY LEMMAS AND OUTLINE OF PROOFS

The aim of this Appendix is to sketch the proofs of the results stated in the previous sections. The full details can be found in a technical report, available from the authors or from http://www.maths.bris.ac.uk/ mapzf/ tvarch/trNLS.pdf.

Before proving these results, we first obtain some results related to weighted sums of tvARCH($p$) processes that we use below.

In what follows, we use $K$ to denote a generic finite positive constant.

**A.1. Properties of tvARCH($p$) processes.** Let us define the following quantity:

$$(25) \qquad \underline{r}(u) = \mathbb{E}\left\{ \frac{\tilde{X}_k^2(u)\tilde{\mathcal{X}}_{k-1}(u)}{\kappa(u_0, \mathcal{X}_{k-1,N})^2} \right\}.$$

LEMMA A.1. *Suppose the conditions in Proposition 3(i) are satisfied, let $\mu(u) = \mathbb{E}\{\tilde{X}_t^2(u)\}$, and let $\mathcal{A}_t(u)$, $\mathcal{D}_t(u)$ and $\underline{r}(u)$ be defined as in (8) and (25), respectively. If $|u_0 - t_0/N| < 1/N$, then we have:*

TABLE 4
*Coverage of 95% Gaussian prediction intervals for our method, using parameter configurations that gave the best results in Table 3*

| Series | CHF | GBP | NOK | NZD | SEK |
|---|---|---|---|---|---|
| Coverage | 99% | 100% | 98% | 79% | 95% |

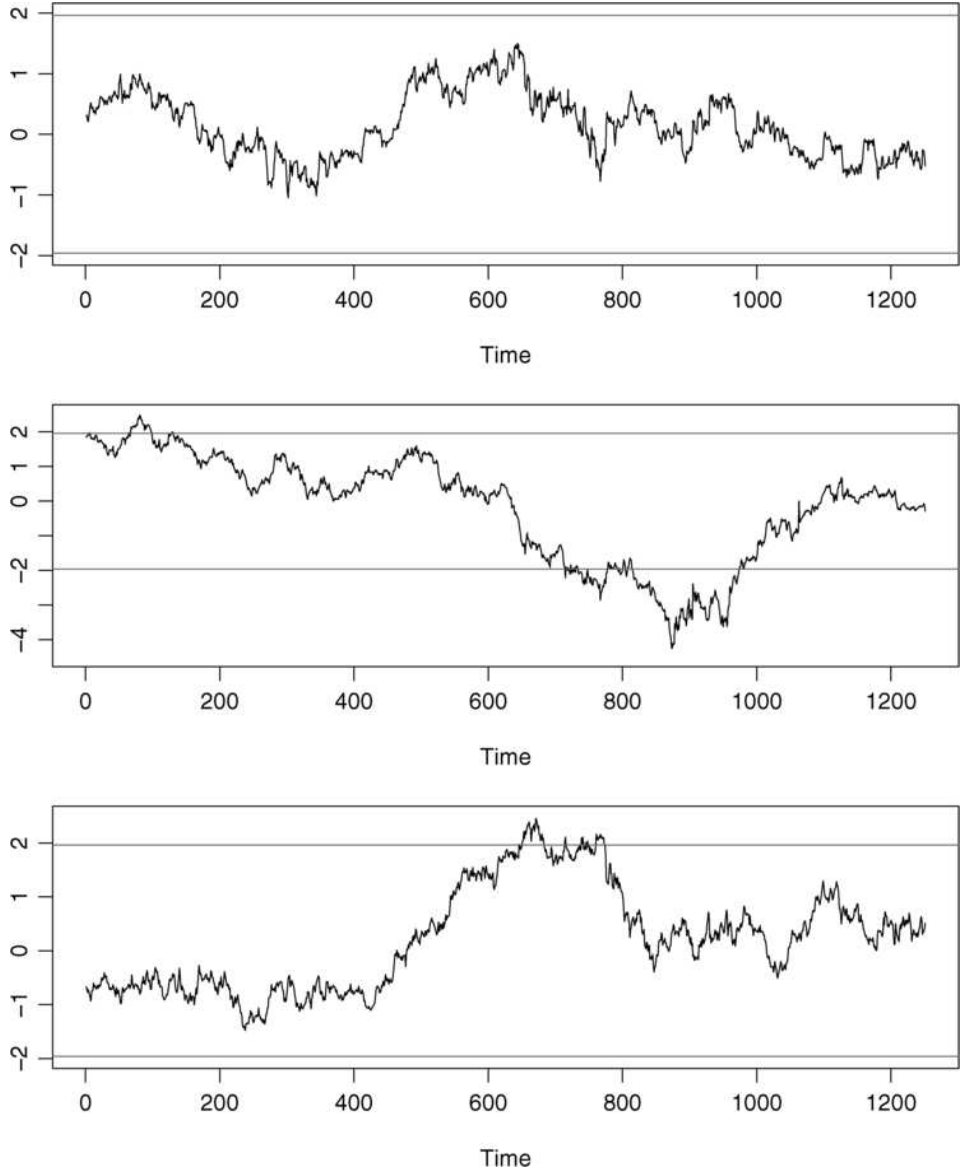

Fig. 4. *From top to bottom: processes $Y_t$ for the GBP, NZD, SEK series. Horizontal lines: symmetric 95% confidence intervals for standard Gaussian variables.*

(i)

$$\sum_{k=1}^{N} \frac{1}{bN} W\left(\frac{t_0 - k}{bN}\right) X_{k,N}^2 \xrightarrow{\mathcal{P}} \mu(u_0).$$
(26)

(ii)

$$\tag{27} \mathcal{R}_{t_0,N} = \sum_{k=p+1}^{N} \frac{1}{bN} W\left(\frac{t_0-k}{bN}\right) \frac{\mathcal{X}_{k-1,N}\mathcal{X}_{k-1,N}^T}{\kappa(u_0,\mathcal{X}_{k-1,N})^2} \xrightarrow{\mathcal{P}} \mathbb{E}[\mathcal{A}_t(u_0)].$$

(iii)

$$\tag{28} \underline{r}_{t_0,N} = \sum_{k=p+1}^{N} \frac{1}{bN} W\left(\frac{t_0-k}{bN}\right) \frac{X_{k,N}^2 \mathcal{X}_{k-1,N}}{\kappa(u_0,\mathcal{X}_{k-1,N})^2} \xrightarrow{\mathcal{P}} \underline{r}(u_0).$$

(iv) *Suppose further that the conditions in Proposition 3(ii) are satisfied, then we have*

$$\tag{29} \sum_{k=p+1}^{N} \frac{1}{bN} W^2\left(\frac{t_0-k}{bN}\right) \frac{\sigma_{k,N}^4 \mathcal{X}_{k-1,N}\mathcal{X}_{k-1,N}^T}{\kappa(u_0,\mathcal{X}_{k-1,N})^4} \xrightarrow{\mathcal{P}} w_2 \mathbb{E}[\mathcal{D}_t(u_0)],$$

*with $b \to 0$, $bN \to \infty$ as $N \to \infty$, where $w_2 = \int_{-1/2}^{1/2} W^2(x)dx$.*

PROOF.  It is straightforward to derive (i), (ii), (iii) and (iv) using Lemma A.5 in Dahlhaus and Subba Rao (2006). We omit the details.  □

To prove Lemma A.3 below we use the following lemma, whose proof is based on mixingale arguments. Suppose $1 \le q < \infty$, and let $\| \cdot \|_q$ denote the $\ell_q$-norm of a vector.

LEMMA A.2.  *Suppose $\{\phi_k \colon k = 1, 2, \ldots\}$ is a stochastic process which satisfies $\mathbb{E}(\phi_k) = 0$ and $\mathbb{E}(\phi_k^q) < \infty$ for some $1 < q \le 2$. Further, let $\mathcal{F}_t = \sigma(\phi_t, \phi_{t-1}, \ldots)$, and suppose that there exists a $\rho \in (0,1)$ such that $\{\mathbb{E}\|\mathbb{E}(\phi_k| \mathcal{F}_{k-j})\|_q^q\}^{1/q} \le K\rho^j$. Then we have*

$$\tag{30} \left\{ \mathbb{E}\left\| \sum_{k=1}^{s} a_k \phi_k \right\|_q^q \right\}^{1/q} \le \frac{K}{1-\rho} \left( \sum_{k=1}^{s} |a_k|^q \right)^{1/q}.$$

In Lemma A.3 below we derive rates of convergence for local sums of a stationary ARCH($p$) process. We use this result to prove the long memory result in Proposition 2.

Let us define the following quantities:

$$\tag{31} \begin{aligned} \mu_1(u,d,h) &= \mathbb{E}\{\tilde{X}_t^2(u)\tilde{X}_{t+h}^2(u+d)\}, \\ c(u,d,h) &= \operatorname{cov}\{\tilde{X}_t^2(u), \tilde{X}_{t+h}^2(u+d)\}, \end{aligned}$$

and set $\mu_1(u, 0, h) = \mu_1(u, h)$ and $c(u, 0, h) = c(u, h)$. Define also the following quantities:

$$\mathcal{S}_{k,bN}(u) = \frac{1}{bN} \sum_{s=kbN}^{(k+1)bN-1} \tilde{X}_s^2(u), \tag{32}$$

$$\mathcal{S}_{k,bN}(u, h, d) = \frac{1}{bN} \sum_{s=kbN}^{(k+1)bN-1} \tilde{X}_s^2(u) \tilde{X}_{s+h}^2(u+d). \tag{33}$$

LEMMA A.3.   *Suppose* $\{\tilde{X}_t(u)\}_t$ *is a stationary* ARCH($p$) *process defined as in* (2) *and suppose the conditions on the parameters* $\{a_j(u)\}_j$ *and the innovations* $\{Z_t\}$ *in Assumption* 1(i), (ii), (iv) *hold. Let* $\mu(u) = \mathbb{E}\{\tilde{X}_t^2(u)\}$, *and let* $\mu_1(u, d, h)$, $\mathcal{S}_{k,bN}(u)$ *and* $\mathcal{S}_{k,bN}(u, h, d)$ *be defined as in* (31), (32) *and* (33) *respectively. Then, we have*

$$\left\{ \mathbb{E} \left\| \sum_{k=p+1}^{N} \frac{1}{bN} W\left(\frac{t-k}{bN}\right) \{\tilde{X}_k^2(u) - \mu(u)\} \right\|_{1+\eta}^{1+\eta} \right\}^{1/(1+\eta)}$$
$$\leq K(bN)^{-(\eta)/1+\eta}. \tag{34}$$

*Further, if* $\{\mathbb{E}(|Z_t|^{2(2+\zeta)})\}^{1/(2+\zeta)} \sup_u \sum_{j=1}^p a_j(u) \leq 1 - \delta$ *for some* $0 < \zeta \leq 2$ *and* $\delta > 0$, *then we have*

$$\{\mathbb{E}\|\mathcal{S}_{k,bN}(u, h, d) - \mu_1(u, d, h)\|_{1+\zeta/2}^{1+\zeta/2}\}^{1/(1+\zeta/2)} \leq K(bN)^{-(\zeta)/2+\zeta}, \tag{35}$$

*where the constant* $K$ *is independent of* $u$ *and* $d$.

PROOF.   We will first prove (34). We use Lemma A.2, with $a_k = W\left(\frac{t-k}{bN}\right)$, $\phi_k = \{\tilde{X}_k^2(u) - \mu(u)\}$ and $q = 1 + \eta$. It can be shown that

$$\{\|\mathbb{E}(\tilde{X}_k^2(u) \mid \mathcal{F}_{k-j}) - \mu(u)\|_{1+\eta}^{1+\eta}\}^{1/(1+\eta)}$$
$$\leq K\rho^j (1 + \{\mathbb{E}\|\tilde{\mathcal{X}}_{k-j}(u)\|_{1+\eta}^{1+\eta}\}^{1/(1+\eta)}),$$

where $\mathcal{F}_t = \sigma(\tilde{X}_t^2(u), \tilde{X}_{t-1}^2(u), \ldots)$. By using the above and that the support of $W\left(\frac{t-k}{bN}\right)$ is proportional to $bN$, we can apply Lemma A.2 to obtain

$$\left\{ \mathbb{E} \left\| \sum_{k=p+1}^{N} \frac{1}{bN} W\left(\frac{t-k}{bN}\right) \{\tilde{X}_k^2(u) - \mu(u)\} \right\|_{1+\eta}^{1+\eta} \right\}^{1/(1+\eta)}$$
$$\leq \frac{1}{bN} \frac{K}{1-\rho} \left( \sum_{k=p+1}^{N} \left| W\left(\frac{t-k}{bN}\right) \right|^{1+\eta} \right)^{1/(1+\eta)} \leq K(bN)^{-\eta/(1+\eta)}.$$

Thus, we have proved (34). The proof of (35) is similar to the proof above but requires the additional (stated) assumption, hence we omit the details. □

**A.2. The covariance structure and the long memory effect of $\mathtt{tvARCH}(p)$ processes.** In this section, we prove results for the covariance structure and the long memory effect of $\mathrm{tvARCH}(p)$ processes.

PROOF OF PROPOSITION 1. It follows easily by making a time-varying Volterra series expansion of the $\mathrm{tvARCH}(p)$ process [see Section 5 in Dahlhaus and Subba Rao (2006)] and using Lemma 2.1 in Giraitis, Kokoszka and Leipus (2000). We omit the details. □

The following lemma is used to prove Proposition 2.

LEMMA A.4. *Suppose $\{X_{t,N}\}_t$ is a $\mathrm{tvARCH}(p)$ which satisfies Assumption 1(i), (ii), (iv), and let $\{\tilde{X}_t(u)\}_t$ be defined as in (2). Let $h := h(N)$ be such that $h/N \to d \in [0,1)$ as $N \to \infty$. Then we have*

$$(36) \qquad \frac{1}{N-h} \sum_{s=1}^{N-h} X_{s,N}^2 \xrightarrow{\mathcal{P}} \int_0^{1-d} \mathbb{E}\{\tilde{X}_t^2(u)\}\, du.$$

*Further, if $\{\mathbb{E}(|Z_t|^{2(2+\zeta)})\}^{1/(2+\zeta)} \sup_u \sum_{j=1}^p a_j(u) \le 1 - \delta$ for some $0 < \zeta \le 2$ and $\delta > 0$, then we have*

$$(37) \qquad \frac{1}{N-h} \sum_{s=1}^{N-h} X_{s,N}^2 X_{s+h,N}^2 \xrightarrow{\mathcal{P}} \int_0^{1-d} \mathbb{E}\{\tilde{X}_t^2(u)\tilde{X}_{t+h}^2(u+d)\}\, du.$$

PROOF. We first prove (37). Let $b := b(N)$ be such that $1/b$ is an integer, $b \to 0$ and $b(N-h) \to \infty$ as $N \to \infty$. We partition the left-hand side of (37) into $1/b$ blocks. Let $k_b = kb[1-d]$, and replace the terms $X_{kb(N-h)+r,N}^2$ and $X_{kb(N-h)+r+h,N}^2$ with $\tilde{X}_{kb(N-h)+r}^2(k_b)$ and $\tilde{X}_{kb(N-h)+r+h}^2(k_b+d)$, respectively. Let $N' = (N-h)$. If $s \in [kbN', (k+1)bN']$ then we replace $X_{s,N}^2$ with $\tilde{X}_s^2(k_b)$ and $X_{s+h,N}^2$ with $\tilde{X}_s^2(k_b+d)$. Recall the notation $\mathcal{S}_{k,bN}(u,h,d)$ and $\mu_1(u,d,h)$ given in (31) and (33), respectively. Now, by using Lemma 1 and that $\frac{kb(N-h)}{N} \le \frac{s}{N} < \frac{(k+1)b(N-h)}{N}$, we have

$$
\begin{aligned}
(38) \qquad & \frac{1}{N-h} \sum_{s=1}^{N-h} X_{s,N}^2 X_{s+h,N}^2 \\
& = b \sum_{k=0}^{1/b-1} \frac{1}{b(N-h)} \sum_{r=0}^{b(N-h)-1} \mathcal{S}_{k,bN}(k_b, h, d) + R_N,
\end{aligned}
$$

where

$$|R_N| \le b \sum_{k=0}^{1/b-1} \frac{1}{bN'} \sum_{r=0}^{bN'-1} \left\{ X_{kbN'+r,N}^2 \left( \frac{1}{N^\beta} V_{kbN'+r+h,N} \right. \right.$$

$$+ \left( 2b + \left| \frac{h}{N} - d \right| \right)^{\beta} W_{kbN'+r+h} \right)$$

$$+ \tilde{X}^2_{kbN'+r+h}(k_b + d) \left( \frac{1}{N^{\beta}} V_{kbN'+r,N} \right.$$

$$\left. + (2b)^{\beta} W_{kbN'+r} \right) \Bigg\}.$$

Now, taking expectations of the above, we have

$$(39) \qquad (\mathbb{E}|R_N|^{1+\zeta/2})^{1/(1+\zeta/2)} \le K \left\{ \left( 2b + b \left| d - \frac{h}{N} \right| \right)^{\beta} + \frac{1}{N^{\beta}} \right\}.$$

Therefore, $R_N \overset{\mathcal{P}}{\to} 0$ as $N \to \infty$. By substituting the integral $\int_0^{1-d} \mu_1(u, d, h)\, du$ with a sum and using (38) and (39), we have

$$\left\{ \mathbb{E} \left\| \frac{1}{N-h} \sum_{s=1}^{N-h} X^2_{s,N} X^2_{s+h,N} - \int_0^{1-d} \mu_1(u, d, h)\, du \right\|_{1+\zeta/2}^{1+\zeta/2} \right\}^{1/(1+\zeta/2)}$$

$$(40) \qquad \le b \sum_{k=0}^{1/b-1} \{ \mathbb{E} \| \mathcal{S}_{k, b(N-h)}(k_b, h, d) - \mu_1(k_b, d, h) \|_{1+\zeta/2}^{1+\zeta/2} \}^{1/(1+\zeta/2)}$$

$$+ \left| b \sum_{k=0}^{1/b-1} \mu_1(k_b, d, h) - \int_0^{1-d} \mu_1(k_b, d, h)\, du \right|$$

$$+ O \left\{ \left( b + kb \left| d - \frac{h}{N} \right| \right)^{\beta} + \frac{1}{N^{\beta}} \right\}.$$

Finally, by substituting the bound (35) and

$$\left| b \sum_{k=0}^{1/b-1} \mu_1(k_b, d, h) - \int_0^{1-d} \mathbb{E}\{ \tilde{X}^2_t(u) \tilde{X}^2_{t+h}(u+d) \}\, du \right| \le kb$$

into (40), we have

$$\left\{ \mathbb{E} \left\| \frac{1}{N-h} \sum_{s=1}^{N-h} X^2_{s,N} X^2_{s+h,N} \right. \right.$$

$$\left. \left. - \int_0^{1-d} \mathbb{E}\{ \tilde{X}^2_t(u) \tilde{X}^2_{t+h}(u+d) \}\, du \right\|_{1+\zeta/2}^{1+\zeta/2} \right\}^{1/(1+\zeta/2)} \to 0,$$

which gives us (37). The proof of (36) is similar and we omit the details. $\square$

PROOF OF PROPOSITION 2.   We first consider the more general case where $h := h(N)$ is such that $h/N \to d \in [0, 1)$ as $N \to \infty$. Then, for fixed $h > 0$, we obtain (5) as special case with $d = 0$.

Let $\mathcal{S}_N(h) = A_N - B_N$, where

$$A_N = \frac{1}{N-h} \sum_{t=1}^{N-h} X_{t,N}^2 X_{t+h,N}^2 \quad \text{and} \quad B_N = (\bar{X}_N)^2.$$

We consider the asymptotic behavior of the terms $A_N$ and $B_N$ separately. By using (36) and (37), we have

$$A_N \xrightarrow{\mathcal{P}} \int_0^{1-d} \mu_1(u, h, d) \, du \quad \text{and} \quad B_N \xrightarrow{\mathcal{P}} \int_0^{1-d} \int_0^{1-d} \mu(u)\mu(v) \, du \, dv.$$

Recall that $\mu(u) = \mathbb{E}\{\tilde{X}_t^2(u)\}$, and that $\mu_1(u, d, h)$ and $c(u, d, h)$ are defined in (31). By using the formula $\mu_1(u, d, h) = c(u, d, h) + \mu(u)\mu(u+d)$, we obtain

$$(41) \quad \mathcal{S}_N(h) \xrightarrow{\mathcal{P}} \int_0^{1-d} \{c(u, d, h) + \mu(u)\mu(u+d)\} \, du - \left\{ \int_0^{1-d} \mu(u) \, du \right\}^2$$

where $h/N \to d$ as $N \to \infty$.

Let us now consider the special case of (41) where $d = 0$. Then, for fixed $h > 0$, we have

$$\mathcal{S}_N(h) \xrightarrow{\mathcal{P}} \int_0^1 c(u, h) \, du + \int_0^1 \int_0^1 \mu^2(u) \, du \, dv - \int_0^1 \int_0^1 \mu(u)\mu(v) \, du \, dv,$$

as $N \to \infty$. This proves (5) and, hence, we have the required result. □

**A.3. Proofs in Section 3.2.** In this section, we prove consistency and asymptotic normality of the weighted kernel-NLS estimator.

PROOF OF PROPOSITION 3(i). By using (27), (28) and Slutsky's theorem, we have

$$\underline{\hat{a}}_{t_0, N} = \{\mathcal{R}_{t_0, N}\}^{-1} \underline{r}_{t_0, N} \xrightarrow{\mathcal{P}} \{\mathbb{E}[\mathcal{A}_t(u_0)]\}^{-1} \underline{r}(u_0).$$

Therefore, to prove that $\underline{\hat{a}}_{t_0, N} \xrightarrow{\mathcal{P}} \underline{a}(u_0)$, we need to show that $\underline{a}(u_0) = \{\mathbb{E}[\mathcal{A}_t(u_0)]\}^{-1} \underline{r}(u_0)$. By using (2) and dividing by $\kappa(u_0, \tilde{\mathcal{X}}_{k-1}(u_0))$, $\{\tilde{X}_{k-1}^2(u_0)\}_k$ satisfies the representation

$$(42) \quad \frac{\tilde{X}_k^2(u_0)}{\kappa(u_0, \tilde{\mathcal{X}}_{k-1}(u_0))} = \frac{\underline{a}^T(u_0)\tilde{\mathcal{X}}_{k-1}(u_0)}{\kappa(u_0, \tilde{\mathcal{X}}_{k-1}(u_0))} + (Z_k^2 - 1)\frac{\tilde{\sigma}_k^2(u_0)}{\kappa(u_0, \tilde{\mathcal{X}}_{k-1}(u_0))}.$$

Finally, multiplying (42) by $\tilde{\mathcal{X}}_{k-1}(u_0)/\kappa(u_0, \mathcal{X}_{k-1, N})$ and taking on both sides expectations, we obtain the desired result.

To prove Proposition 3(ii), we use the same methodology given in the proof of Theorem 3 in Dahlhaus and Subba Rao (2006), hence we omit the details. □

**A.4. Proofs in Section 3.4.** In this section, we prove consistency and asymptotic normality of the two-stage kernel-NLS estimator. To prove these asymptotic properties, we need the following two lemmas.

LEMMA A.5. *Suppose* $\{X_{t,N}\}_t$ *is a* tvARCH($p$) *process which satisfies Assumption* 1(i), (ii), (iv), *let* $\mu(u) = \mathbb{E}\{\tilde{X}_t^2(u)\}$, *and let* $\hat{\mu}_{t_0,N}$ *be defined as in* (15). *If* $|u_0 - t_0/N| < 1/N$, *then, for* $0 \leq i, j \leq p$, *we have*

$$(43) \quad \sum_{k=p+1}^{N} \frac{1}{bN} W\left(\frac{t_0 - k}{bN}\right) \frac{X_{k-i,N}^2 X_{k-j,N}^2}{(\hat{\mu}_{t_0,N} + \mathcal{S}_{k-1,N})^2} \xrightarrow{\mathcal{P}} \mathbb{E}\left(\frac{\tilde{X}_{k-i}^2(u_0)\tilde{X}_{k-j}^2(u_0)}{(\mu(u_0) + \mathcal{S}_{k-1}(u_0))^2}\right),$$

*with* $b \to 0$, $bN \to \infty$ *as* $N \to \infty$.

PROOF. To prove the result we use techniques similar to those in Bose and Mukherjee (2003). By using the inequality $|1/x^2 - 1/y^2| \leq 2|x-y|\{(1/x)[1+x/y]\}^3$, for $x, y > 0$, we bound the difference

$$(44) \quad \begin{aligned} &\left| \frac{X_{k-i,N}^2 X_{k-j,N}^2}{(\hat{\mu}_{t_0,N} + \mathcal{S}_{k-1,N})^2} - \frac{X_{k-i,N}^2 X_{k-j,N}^2}{(\mu(u_0) + \mathcal{S}_{k-1,N})^2} \right| \\ &\leq 2 X_{k-i,N}^2 X_{k-j,N}^2 |\hat{\mu}_{t_0,N} - \mu(u_0)| \\ &\quad \times \left| \frac{1}{\mu(u_0) + \mathcal{S}_{k-1,N}} \left(1 + \frac{\hat{\mu}_{t_0,N} + \mathcal{S}_{k-1,N}}{\mu(u_0) + \mathcal{S}_{k-1,N}}\right) \right|^3 \\ &\leq 2|\hat{\mu}_{t_0,N} - \mu(u_0)| \left(1 + \frac{|\mu(u_0)|}{|\hat{\mu}_{t_0,N}|}\right) \frac{X_{k-i,N}^2 X_{k-j,N}^2}{(\mu(u_0) + \mathcal{S}_{k-1,N})^3}. \end{aligned}$$

Let us now define the following quantities:

$$\Gamma(u_0) = \mathbb{E}\left(\frac{\tilde{X}_{k-i}^2(u_0)\tilde{X}_{k-j}^2(u_0)}{\{\mu(u_0) + \mathcal{S}_{k-1}(u_0)\}^2}\right),$$

$$A_{t_0,N} = \sum_{k=p+1}^{N} \frac{1}{bN} W\left(\frac{t_0 - k}{bN}\right) \frac{X_{k-i,N}^2 X_{k-j,N}^2}{\{\hat{\mu}_{t_0,N} + \mathcal{S}_{k-1,N}\}^2},$$

$$C_{t_0,N}(u_0) = \left| \sum_{k=p+1}^{N} \frac{1}{bN} W\left(\frac{t_0 - k}{bN}\right) \frac{X_{k-i,N}^2 X_{k-j,N}^2}{\{\mu(u_0) + \mathcal{S}_{k-1,N}\}^2} - \Gamma(u_0) \right|.$$

Then, by using the bound (44), we have

$$\begin{aligned} &\left| A_{t_0,N} - \Gamma(u_0) \right| \\ &\leq \left| \sum_{k=p+1}^{N} \frac{1}{bN} W\left(\frac{t_0 - k}{bN}\right) \frac{X_{k-i,N}^2 X_{k-j,N}^2}{\{\mu(u_0) + \mathcal{S}_{k-1,N}\}^2} - \Gamma(u_0) \right| \end{aligned}$$

$$(45) \qquad + 2|\hat{\mu}_{t_0,N} - \mu(u_0)|\left(1 + \frac{|\mu(u_0)|}{|\hat{\mu}_{t_0,N}|}\right)$$

$$\times \sum_{k=p+1}^{N} \frac{1}{bN}\left|W\left(\frac{t_0-k}{bN}\right)\right|\frac{X_{k-i,N}^2 X_{k-j,N}^2}{\{\mu(u_0) + \mathcal{S}_{k-1,N}\}^3}$$

$$\leq C_{t_0,N}(u_0) + 2|\hat{\mu}_{t_0,N} - \mu(u_0)|\left(1 + \frac{|\mu(u_0)|}{|\hat{\mu}_{t_0,N}|}\right)\sum_{k=p+1}^{N} \frac{1}{bN}\left|W\left(\frac{t_0-k}{bN}\right)\right|.$$

Since $\hat{\mu}_{t_0,N} \xrightarrow{\mathcal{P}} \mu(u_0)$, by using Slutsky's lemma we have

$$|\hat{\mu}_{t_0,N} - \mu(u_0)|\left(1 + \frac{|\mu(u_0)|}{|\hat{\mu}_{t_0,N}|}\right) \xrightarrow{\mathcal{P}} 0.$$

Furthermore, by using (27) we have $C_{t_0,N}(u_0) \xrightarrow{\mathcal{P}} 0$. Altogether this gives $|A_{t_0,N} - \Gamma(u_0)| \xrightarrow{\mathcal{P}} 0$, and the desired result follows. $\square$

To show asymptotic normality, we need to define the following least-squares criteria:

$$(46) \qquad \tilde{\mathcal{L}}_{t_0,N}(\underline{\alpha}) = \sum_{k=p+1}^{N} \frac{1}{bN}W\left(\frac{t_0-k}{bN}\right)\tilde{h}_{t_0,N}(X_{k,N}, \underline{X}_{k-1,N}, \underline{\alpha}),$$

$$(47) \qquad \tilde{\mathcal{L}}_{t_0}(u, \underline{\alpha}) = \sum_{k=p+1}^{N} \frac{1}{bN}W\left(\frac{t_0-k}{bN}\right)\tilde{h}_{t_0,N}(\tilde{X}_k(u), \underline{\tilde{X}}_{k-1}(u), \underline{\alpha}),$$

$$(48) \qquad \mathcal{L}_{t_0}^{(\mu)}(u, \underline{\alpha}) = \sum_{k=p+1}^{N} \frac{1}{bN}W\left(\frac{t_0-k}{bN}\right)\tilde{h}(u, \tilde{X}_k(u), \underline{\tilde{X}}_{k-1}(u), \underline{\alpha}),$$

where $\tilde{h}_{t_0,N}(y_0, \underline{y}, \underline{\alpha}) = (\hat{\mu}_{t_0,N} + \sum_{j=1}^{p} y_j^2)^{-2}\{y_p^2 - \alpha_0 - \sum_{j=1}^{p} \alpha_j y_j^2\}^2$ and $\tilde{h}(u, y_0, \underline{y}, \underline{\alpha}) = (\mu(u) + \sum_{j=1}^{p} y_j^2)^{-2}\{y_0^2 - \alpha_0 - \sum_{j=1}^{p} \alpha_j y_j^2\}^2$. We note that $\underline{\tilde{a}}_{t_0,N} = \arg\min_{\underline{a}} \tilde{\mathcal{L}}_{t_0,N}(\underline{a})$. Asymptotic normality of $\sqrt{bN}\nabla \mathcal{L}_{t_0}^{(\mu)}(u_0, \underline{a}(u_0))$ can easily be established by verifying the conditions of the martingale central limit theorem. However, the same theorem cannot be used to show the asymptotic normality of $\sqrt{bN}\nabla \tilde{\mathcal{L}}_{t_0}(u_0, \underline{a}(u_0))$, since $\tilde{\mathcal{L}}_{t_0}(u_0, \underline{a}(u_0))$ is not a sum of martingale differences. In Lemma A.6 below we overcome this problem by showing that $\sqrt{bN}(\tilde{\mathcal{L}}_{t_0}(u_0, \underline{a}(u_0)) - \mathcal{L}_{t_0}^{(\mu)}(u_0, \underline{a}(u_0))) \xrightarrow{\mathcal{P}} 0$, which allows us to replace $\tilde{\mathcal{L}}_{t_0}(u_0, \underline{a}(u_0))$ with $\mathcal{L}_{t_0}^{(\mu)}(u_0, \underline{a}(u_0))$.

LEMMA A.6. *Suppose* $\{X_{t,N}\}_t$ *is a* tvARCH($p$) *process which satisfies Assumption* 1(i), (ii), (iii), (iv). *Let* $\mu(u) = \mathbb{E}\{\tilde{X}_t^2(u)\}$, *and* $\hat{\mu}_{t_0,N}$, $\tilde{\mathcal{L}}_{t_0}(u, \underline{\alpha})$ *and $\mathcal{L}_{t_0}^{(\mu)}(u, \underline{\alpha})$ be defined as in* (15), (47) *and* (48), *respectively. If $|u_0 - t_0| < 1/N$, then we have*

$$(49) \qquad |\hat{\mu}_{t_0, N} - \mu(u_0)| = O_p(b^\beta + (bN)^{-\eta/(1+\eta)}),$$

*and*

$$(50) \qquad \sqrt{bN}[\nabla\tilde{\mathcal{L}}_{t_0}(u_0, \underline{a}(u_0)) - \nabla\mathcal{L}_{t_0}^{(\mu)}(u_0, \underline{a}(u_0))] = o_p(1),$$

*with $b \to 0$, $bN \to \infty$ as $N \to \infty$.*

PROOF OF PROPOSITION 4.   (i) It is straightforward to show consistency using (17) and Lemma A.5.

(ii) Define $\tilde{\mathcal{B}}_{t_0, N}(\underline{\alpha}) = \tilde{\mathcal{L}}_{t_0, N}(\underline{\alpha}) - \tilde{\mathcal{L}}_{t_0}(u, \underline{\alpha})$. To prove that $\nabla\tilde{\mathcal{B}}_{t_0, N}(\underline{a}(u_0)) = O_p(b^\beta)$, we use the same arguments given in Theorem 3 in Dahlhaus and Subba Rao (2006), hence we omit the details. To prove (17), we use that $\nabla\tilde{\mathcal{L}}_{t_0, N}(\tilde{\underline{a}}_{t_0, N}) = 0$, and expanding $\nabla\tilde{\mathcal{L}}_{t_0, N}(\tilde{\underline{a}}_{t_0, N})$ about $\underline{a}(u_0)$, we have

$$\nabla^2\tilde{\mathcal{L}}_{t_0, N}(\tilde{\underline{a}}_{t_0, N} - \underline{a}(u_0))$$
$$= \mathcal{L}_{t_0}^{(\mu)}(u_0, \underline{a}(u_0)) + \{\nabla\tilde{\mathcal{L}}_{t_0}(u_0, \underline{a}(u_0)) - \nabla\mathcal{L}_{t_0}^{(\mu)}(u_0, \underline{a}(u_0))\}$$
$$- \nabla\tilde{\mathcal{B}}_{t_0, N}(\underline{a}(u_0)).$$

By using Lemma A.5, we easily see that $\nabla^2\tilde{\mathcal{L}}_{t_0, N}(\underline{a}(u_0)) \xrightarrow{\mathcal{P}} 2\mathbb{E}[\mathcal{A}_t^{(\mu)}(u_0)]$ and, using (50), we have

$$(\tilde{\underline{a}}_{t_0, N} - \underline{a}(u_0))$$
$$= -\frac{1}{2}\{\nabla\mathcal{L}_{t_0}^{(\mu)}(u_0, \underline{a}(u_0)) + \nabla\tilde{\mathcal{B}}_{t_0, N}(\underline{a}(u_0))\}\{\mathbb{E}[\mathcal{A}_t^{(\mu)}(u_0)]\}^{-1} + o_p\left(\frac{1}{bN}\right).$$

Finally, by using the martingale central limit theorem [see, e.g., Hall and Heyde (1980), Theorem 3.2], we obtain (17).   □

**A.5. Proofs in Section 5.**   In this section, we prove the results in Section 5. Some of the results in this section have been inspired by corresponding results in the residual bootstrap for linear processes literature [cf. Franke and Kreiss (1992)]. However, the proofs are technically very different, because the tvARCH($p$) process is a nonlinear, nonstationary process, and the normalization of the two-stage kernel-NLS estimator with random weights.

In order to show that the distribution of the bootstrap sample $\hat{\underline{a}}_{t_0, N}^+ - \bar{\underline{a}}_{t_0, N}$ asymptotically coincides with the asymptotic distribution of $\tilde{\underline{a}}_{t_0, N} - \underline{a}(u_0)$, we will show convergence of the distributions under the Mallows distance. The Mallows distance between the distribution $H$ and $G$ is defined as

$$d_2(H, G) = \inf_{X \sim H, Y \sim G}\{\mathbb{E}(X - Y)^2\}^{1/2}.$$

Roughly speaking, if $d_2(F_n, G_n) \to 0$, then the limiting distributions of $F_n$ and $G_n$ are the same [Bickel and Freedman (1981)]. Following Franke and Kreiss (1992), to reduce notation, we let $d_2(X, Y) = d_2(H, G)$, where the random variables $X$ and $Y$ have measures $H$ and $G$, respectively.

We also require the following definitions. Let

$$\tilde{\mathcal{R}}_N(u_0) = \sum_{k=p+1}^{N} \frac{1}{bN} W\left(\frac{t_0 - k}{bN}\right) \frac{\tilde{\mathcal{X}}_{k-1}(u_0) \tilde{\mathcal{X}}_{k-1}(u_0)^T}{(\hat{\mu}_{t_0, N} + \sum_{j=1}^{k} X_{k-j}^2(u_0))},$$

$$\tilde{\underline{r}}_N(u_0) = \sum_{k=p+1}^{N} \frac{1}{bN} W\left(\frac{t_0 - k}{bN}\right) \frac{\tilde{X}_k^2(u_0) \tilde{\mathcal{X}}_{k-1}(u_0)^T}{(\hat{\mu}_{t_0, N} + \sum_{j=1}^{k} \tilde{X}_{k-j}^2(u_0))^2}.$$

PROPOSITION 6. *Suppose Assumption 1 holds, and suppose either* $\inf_j a_j(u_0) > 0$ *or* $\mathbb{E}(Z_t^4)^{1/2} \sup_u[\sum_{j=1}^{p} a_j(u)] \le 1 - \delta$ *[which implies* $\sup_k \mathbb{E}(X_{k,N}^4) < \infty$]. *Let* $F$ *be the distribution function of* $Z_t^2$. *Then we have*

$$(51) \qquad\qquad d_2(\hat{F}_{t_0, N}, F) \xrightarrow{\mathcal{P}} 0.$$

*Furthermore, if we suppose* $b^\beta \sqrt{bN} \to 0$, *then we have*

$$(52) \quad d_2(\sqrt{bN}(\underline{r}_{t_0, N}^+ - \mathcal{R}_{t_0, N}^+ \underline{a}_{0, N}), \sqrt{bN}(\tilde{\underline{r}}_N(u_0) - \tilde{\mathcal{R}}_N(u_0)\underline{a}(u_0))) \xrightarrow{\mathcal{P}} 0,$$

*and*

$$(53) \qquad\qquad \mathcal{R}_{t_0, N}^+ \xrightarrow{\mathcal{P}} \mathbb{E}\{A_t^{(\mu)}(u_0)\},$$

*with* $b \to 0$, $bN \to \infty$ *as* $N \to \infty$.

We prove each part of the proposition below.

PROOF OF (51). To prove this result, we define the empirical distribution function of the true residuals, that is,

$$F_{t_0, N}(x) = \frac{1}{2bN} \sum_{k=t_0 - bN}^{t_0 + bN - 1} \mathbb{I}_{(-\infty, x]}(Z_t^2),$$

noting that $\hat{Z}_t^2$ is an estimator of $Z_t^2$. [It is worth pointing out that in a different context, the empirical distribution of the estimated residuals of a stationary ARCH($p$) process was considered in Horváth, Kokoszka and Teyssiére (2001).] We first observe that since $d_2$ is a distance it satisfies the triangle inequality $d_2(\hat{F}_{t_0, N}, F) \le d_2(\hat{F}_{t_0, N}, F_{t_0, N}) + d_2(F_{t_0, N}, F)$. By using Lemma 8.4 in Bickel and Freedman (1981), it can be shown that $d_2(F_{t_0, N}, F) \xrightarrow{\mathcal{P}} 0$. Therefore, to prove (51), we need only show that $d_2(\hat{F}_{t_0, N}, F_{t_0, N}) \xrightarrow{\mathcal{P}} 0$.

By definition of $d_2$ and the measures $\hat{F}_{t_0,N}$ and $F_{t_0,N}$, we have

$$d_2(\hat{F}_{t_0,N}, F_{t_0,N})^2 = \inf_{Z_t^{+2} \in \hat{F}_{t_0,N}, Z_t^2 \in F_{t_0,N}} \mathbb{E}(Z_t^{+2} - Z_t^2)^2,$$

where the infimum is taken over all joint distributions on $(Z_t^{+2}, Z_t^2)$ which have marginals $\hat{F}_{t_0,N}$ and $F_{t_0,N}$. Let us suppose $P(J = i) = (i + bN)/2bN$, for $i \in \{-bN, \ldots, bN - 1\}$, and define $Z_t^{+2} = \hat{Z}_J^2$ and $Z_t^2 = Z_J^2$. Then, since $(\hat{Z}_J^2, Z_J^2)$ both have marginals $\hat{F}_{t_0,N}$ and $F_{t_0,N}$, respectively, we have

$$d_2(\hat{F}_{t_0,N}, F_{t_0,N})^2 \leq \mathbb{E}(\hat{Z}_J^2 - Z_J^2)^2 = \frac{1}{2bN} \sum_{k=t_0-bN}^{t_0+bN-1} (\hat{Z}_k^2 - Z_k^2)^2$$

$$\leq \frac{1}{2bN} \sum_{k=t_0-bN}^{t_0+bN-1} \left( Z_k^2 - \hat{Z}_k^2 + \frac{1}{2bN} \sum_{k=t_0-bN}^{t_0+bN} \hat{Z}_k^2 - 1 \right)^2.$$

By adding and subtracting $\frac{1}{2bN} \sum_{k=t_0-bN}^{t_0+bN} Z_k^2$, and using that $\frac{1}{2bN} \sum_{k=t_0-bN}^{t_0+bN} (Z_k^2 - 1)^2 \xrightarrow{\mathcal{P}} 0$, we have

$$d_2(\hat{F}_{t_0,N}, F_{t_0,N})^2$$

$$\leq \frac{1}{2bN} \sum_{k=t_0-bN}^{t_0+bN-1} \left( Z_k^2 - \hat{Z}_k^2 + \frac{1}{2bN} \sum_{k=t_0-bN}^{t_0+bN-1} (\hat{Z}_k^2 - Z_k^2) \right.$$

$$\left. + \frac{1}{2bN} \sum_{k=t_0-bN}^{t_0+bN-1} (Z_k^2 - 1) \right)^2$$

$$\leq \frac{K}{bN} \sum_{k=t_0-bN}^{t_0+bN-1} (\hat{Z}_k^2 - Z_k^2)^2 + \frac{K}{bN} \sum_{k=t_0-bN}^{t_0+bN-1} (Z_k^2 - 1)^2$$

$$\leq \frac{K}{bN} \sum_{k=t_0-bN}^{t_0+bN-1} \left[ \frac{Z_k^2}{\hat{\sigma}_{k,N}^2} \left\{ a_0(k/N) - \tilde{a}_{t_0,N}(0) \right. \right.$$

$$\left. \left. + \sum_{j=1}^{p} [a_j(k/N) - \tilde{a}_{t_0,N}(j)] X_{k-j,N}^2 \right\} \right]^2 + o_p(1),$$

where $\hat{\sigma}_{k,N}^2 = \tilde{a}_{t_0,N}(0) + \sum_{j=1}^{p} \tilde{a}_{t_0,N}(j) X_{k-j,N}^2$. Now by bounding the above in two different ways we obtain

$$d_2(\hat{F}_{t_0,N}, F_{t_0,N})^2$$

$$\leq \min \left\{ \begin{array}{l} K\left(\sum_{j=0}^{p} \dfrac{|a_j(k/N) - \tilde{a}_{t_0,N}(j)|}{|\tilde{a}_{t_0,N}(j)|}\right)^2 \left(\dfrac{1}{2bN}\sum_{k=t_0-bN}^{t_0+bN/2} Z_k^4\right) \\[2.5ex] K\left(\sum_{j=0}^{p} \dfrac{|a_j(k/N) - \tilde{a}_{t_0,N}(j)|}{|\tilde{a}_{t_0,N}(0)|}\right)^2 \left(\dfrac{1}{2bN}\sum_{k=t_0-bN}^{t_0+bN/2} Z_k^4 X_{k-j,N}^4\right) \end{array} \right\}$$

$$+ \, o_p(1).$$

To show (51), we need to use the bounds above noting that the bound we use depends on the conditions we have placed on the parameters $\{a_j(\cdot)\}$. By using $|a_j(u) - a_j(v)| \leq K|u-v|^\beta$ and $k \in [t_0 - bN, t_0 + bN - 1]$, we have

$$|\underline{a}(k/N) - \underline{\tilde{a}}_{t_0,N}| \leq |\underline{a}(k/N) - \underline{a}(u_0)| + |\underline{a}(u_0) - \underline{\tilde{a}}_{t_0,N}|$$

$$\leq Kb^\beta + |\underline{a}(u_0) - \underline{\tilde{a}}_{t_0,N}|.$$

Since $|\underline{\tilde{a}}_{t_0,N} - \underline{a}(u_0)| \overset{\mathcal{P}}{\to} 0$, by using the above, it is straightforward to show

$$Kp\left(\frac{b^\beta + |\underline{a}(u_0) - \underline{\tilde{a}}_{t_0,N}|}{\min_j |\tilde{a}_{t_0,N}(j)|}\right)^2 \left(\frac{1}{2bN}\sum_{k=t_0-bN}^{t_0+bN-1} Z_k^4\right) \overset{\mathcal{P}}{\to} 0$$

$$\text{if } \inf_j a_j(u) > 0,$$

$$Kp\left(\frac{b^\beta + |\underline{a}(u_0) - \underline{\tilde{a}}_{t_0,N}|}{|\tilde{a}_{t_0,N}(0)|}\right)^2 \left(\frac{1}{2bN}\sum_{k=t_0-bN}^{t_0+bN-1} Z_k^4 X_{k-j,N}^4\right) \overset{\mathcal{P}}{\to} 0$$

$$\text{if } \sup_{k,N} \mathbb{E}(X_{k,N}^4) < \infty,$$

with $b \to 0$, $bN \to \infty$ as $N \to \infty$. Therefore, under the stated assumptions, and by using the above convergence in probability, we have that $d_2(\hat{F}_{t_0,N}, F_{t_0,N}) \overset{\mathcal{P}}{\to} 0$. Altogether this means that $d_2(\hat{F}_{t_0,N}, F) \overset{\mathcal{P}}{\to} 0$, with $b \to 0$, $bN \to \infty$ as $N \to \infty$, thus we obtain the result. $\square$

It follows from the above [Bickel and Freedman (1981), Lemma 8.3] that

(54)
$$\mathbb{E}(Z_t^{+2}) \overset{\mathcal{P}}{\to} \mathbb{E}(Z_t), \qquad \mathbb{E}(Z_t^{+4}) \overset{\mathcal{P}}{\to} \mathbb{E}(Z_t^4) \quad \text{and}$$

$$\inf \mathbb{E}(Z_t^{+2} - Z_t^2)^2 \overset{\mathcal{P}}{\to} 0,$$

where the infimum is taken over all joint distributions on $(Z_t^{+2}, Z_t^2)$ which have marginals $\hat{F}_{t_0,N}$ and $F_{t_0,N}$. We use these limits to prove the results below.

To prove Proposition 6 we require the following definitions:

$$\tilde{X}_t^{+2}(u_0) = \tilde{\sigma}_t^{+2}(u_0) Z_t^{+2}, \qquad \tilde{\sigma}_t^{+2}(u_0) = a_0(u_0) + \sum_{j=1}^{p} a_j(u_0) \tilde{X}_{t-j}^{+2}(u_0),$$

and Lemma A.7 below. We note that $\tilde{X}_t^{+2}(u_0)$ is very similar to $\tilde{X}_t^{+2}(u_0)$, but the estimated parameters $\underline{\bar{a}}_{t_0,N}$ have been replaced by the true parameters $\underline{a}(u_0)$.

In the lemma below we show that for $t \in [t_0 - bN/2, t_0 + bN/2 - 1]$, the distributions of $X_t^{+2}(u_0)$ and $X_t(u_0)$ are sufficiently close and the difference is uniformly bounded over $t$.

LEMMA A.7. *Suppose assumptions in Proposition 6 hold, then we have*

$$(55) \qquad \mathbb{E}|X_t^{+2}(u_0) - \tilde{X}_t^{+2}(u_0)| \le C|\bar{a}_{t_0,N} - \underline{a}(u_0)| \sum_{k=1}^{\infty} k^2 (1-\delta)^k \xrightarrow{\mathcal{P}} 0,$$

*where $b \to 0$, $bN \to \infty$ as $N \to \infty$ and where the expectation is conditioned on $\{X_{k,N}\}$. Furthermore for $t_0 + bN/2 \le t \le t_0 + bN/2$ we have*

$$
\begin{aligned}
&\inf \mathbb{E}|\tilde{X}_t^{+2}(u_0) - \tilde{X}_t^2(u_0)| \\
(56) \quad &\le C \sum_{k=1}^{\infty} (1 + (\mathbb{E}|Z_t^{+2}|)^{k+bN/(2p)})(1-\delta)^k \\
&\quad + C \inf \mathbb{E}|Z_1^{+2} - Z_t^2| \sum_{k=1}^{\infty} \{1 + \mathbb{E}(Z_1^{+2}) + \cdots + [\mathbb{E}(Z_1^{+2})]^{k-1}\}(1-\delta)^k \\
&\quad + o_p(1) \xrightarrow{\mathcal{P}} 0,
\end{aligned}
$$

*where $b \to 0$, $bN \to \infty$ as $N \to \infty$. The expectation is with respect to the measure on all independent pairs $\{(Z_t^{+2}, Z_t^2)\}_t$, and the infimum is taken over all joint distributions on $(Z_t^{+2}, Z_t^2)$ which have marginals $\hat{F}_{t_0,N}$ and $F_{t_0,N}$, respectively.*

PROOF. It can be shown that the stationary ARCH($\infty$) process has a solution which can be written in terms of a Volterra series [Giraitis, Kokoszka and Leipus (2000)]. Define for all $j > p$, $a_j(u_0) = 0$ and $\bar{a}_{t_0,N}(j) = 0$. Then, by following Giraitis, Kokoszka and Leipus (2000), $X_t^{+2}(u_0), \tilde{X}_t^{+2}(u_0)$ and $\tilde{X}_t^2(u_0)$ have the solutions

$$X_t^{+2}(u_0) = \sum_{k=0}^{N} \sum_{j_k < \cdots j_0 : j_0 = t} \left\{ \prod_{s=0}^{k} \bar{a}_{t_0,N}(j_s - j_{s+1}) \right\} \prod_{s=1}^{k} Z_{j_s}^{+2},$$

$$\tilde{X}_t^{+2}(u_0) = \sum_{k=0}^{N} \sum_{j_k < \cdots j_0 : j_0 = t} \left\{ \prod_{s=0}^{k} a_{j_s - j_{s+1}}(u_0) \right\} \prod_{s=1}^{k} Z_{j_s}^{+2},$$

$$\tilde{X}_t^2(u_0) = \sum_{k=0}^{\infty} \sum_{j_k < \cdots j_0 : j_0 = t} \left\{ \prod_{s=0}^{k} a_{j_s - j_{s+1}}(u_0) \right\} \prod_{s=1}^{k} Z_{j_s}^2,$$

respectively. We first consider

$$\mathbb{E}|X_t^{+2}(u_0) - \tilde{X}_t^{+2}(u_0)|$$

$$= \mathbb{E}\left|\sum_{k=0}^{N}\sum_{j_k < \cdots < j_0; j_0 = t}\left\{\bar{a}_{t_0,N}(j_s - j_{s+1}) - \prod_{s=0}^{k}a_{j_s - j_{s+1}}(u_0)\right\}\prod_{s=1}^{k}Z_{j_s}^{+2}\right|.$$

Now, by repeatedly taking differences, and using that $\sup_u \sum_{j=1}^{p}a_j(u) \leq 1 - \delta$, $\sum_{j=1}^{p}\bar{a}_{t_0,N}(j) \leq 1 - \delta$ and $\|\bar{\underline{a}}_{t_0,N} - \underline{a}(u_0)\|_2 \xrightarrow{\mathcal{P}} 0$, we obtain (55).

To prove (56), we first note that expectation is taken with respect to the joint measure on the independent pairs $\{(Z_t^{+2}, Z_t^2)\}_t$. Using the Volterra expansions above, we have

$$\mathbb{E}|\tilde{X}_t^{+2}(u_0) - \tilde{X}_t^2(u_0)|$$

$$(57)\qquad \leq \sum_{k=0}^{\infty}\sum_{1\leq j_1,\ldots,j_k \leq p}\prod_{s=0}^{k}a_{j_s}(u_0)\mathbb{E}\left|\prod_{s=1}^{k}Z_{t-\sum_{i=1}^{s}j_i}^{+2} - \prod_{s=1}^{k}Z_{t-\sum_{i=1}^{s}j_i}^2\right|$$

$$+ o_p(1).$$

We see from (54), if $t_0 - bN \leq k \leq t_0 + bN - 1$ and by setting $Z_k^{+2} = \hat{Z}_k$, we have $\mathbb{E}|Z_k^{+2} - Z_k^2| \xrightarrow{\mathcal{P}} 0$. Therefore, for all $t \in [t_0 - bN/2, t_0 + bN/2 - 1]$, and $t - bN/2 \leq i \leq t$, we will show that $\inf \mathbb{E}|\prod_{s=1}^{k}Z_{t-\sum_{i=1}^{s}j_i}^{+2} - \prod_{s=1}^{k}Z_{t-\sum_{i=1}^{s}j_i}^2| \xrightarrow{\mathcal{P}} 0$. This allows us to obtain a uniform rate of convergence for $\mathbb{E}|\tilde{X}_t^{+2}(u_0) - \tilde{X}_t^2(u_0)|$ for all $t_0 - bN/2 \leq k \leq t_0 + bN/2 - 1$. To obtain this rate, we partition the inner sum above into two sums, where $\sum_{i=1}^{k}j_s \leq bN/2$ and $\sum_{i=1}^{k}j_s > bN/2$. We further note that since for all $i$, $1 \leq j_i \leq p$, if $\sum_{j=1}^{k}j_s > bN/2$, then this implies $k > bN/(2p)$. Altogether this gives

$$(58)\qquad \mathbb{E}|\tilde{X}_t^{+2}(u_0) - \tilde{X}_t^2(u_0)| \leq I + II + o_p(1),$$

where

$$I = a_0(u_0)\sum_{k=0}^{\infty}\sum_{\sum_{i=1}^{k}j_s \leq bN/2}\prod_{s=1}^{k}a_{j_s}(u_0)\mathbb{E}\left|\prod_{s=1}^{k}Z_{t-\sum_{i=1}^{s}j_i}^{+2} - \prod_{s=1}^{k}Z_{t-\sum_{i=1}^{s}j_i}^2\right|$$

and

$$II = \sum_{k>bN/(2p)}\sum_{1\leq j_1,\ldots,j_k \leq p}\prod_{s=0}^{k}a_{j_s}(u_0)\{(\mathbb{E}|Z_t^2|)^k + (\mathbb{E}|Z_t^{+2}|)^k\}.$$

We now study $I$ and consider, in particular, the difference $\mathbb{E}|\prod_{s=1}^{k}Z_{j_s}^{+2} - \prod_{s=1}^{k}Z_{j_s}^2|$. By repeatedly taking differences, we have

$$\mathbb{E}\left|\prod_{s=1}^{k}Z_{j_s}^{+2} - \prod_{s=1}^{k}Z_{j_s}^2\right| \leq \mathbb{E}|Z_{j_s}^{+2} - Z_{j_s}^2|\{1 + \mathbb{E}(Z_{j_s}^{+2}) + \cdots + [\mathbb{E}(Z_{j_s}^2)]^{k-1}\}.$$

Substituting the above into $I$, taking the infimum over all joint measures on $(Z_t^{+2}, Z_t^2)$, and using $\sup_u \sum_{j=1}^p a_j(u) \leq 1 - \delta$, we obtain

$$a_0(u_0) \sum_{k=0}^{\infty} \sum_{\sum_{i=1}^k j_s \leq bN/2} \prod_{s=0}^k a_{j_s}(u_0) \inf \mathbb{E} \left| \prod_{s=1}^k Z_{t-\sum_{i=1}^s j_i}^{+2} - \prod_{s=1}^k Z_{t-\sum_{i=1}^s j_i}^2 \right|$$

$$(59) \qquad \leq C \{ \inf \mathbb{E} | Z_t^{+2} - Z_t^2 | \}$$

$$\times \sum_{k=1}^N \{ 1 + \mathbb{E}(Z_1^{+2}) + \cdots + [\mathbb{E}(Z_1^{+2})]^k \} (1 - \delta)^k + o_p(1).$$

We note that in the above we have extended the sum beyond $\sum_{i=1}^k j_s \leq bN/2$ to make the summands easier to handle. Our aim is to show that the right-hand side of (59) converges in probability to 0. For any $\varepsilon > 0$, define $B_N^\varepsilon := \{ \mathbb{E} | Z_1^{+2} | > 1 + \varepsilon \}$. By (54), we have $P(B_N^\varepsilon) \to 0$ as $N \to \infty$. Denote further

$$A_n^\varepsilon := \left\{ C \{ \inf \mathbb{E} | Z_t^{+2} - Z_t^2 | \} \sum_{k=1}^N \{ 1 + \mathbb{E}(Z_1^{+2}) + \cdots + [\mathbb{E}(Z_1^{+2})]^k \} (1 - \delta)^k > \varepsilon \right\}.$$

For $\varepsilon_1 < \delta/(1 - \delta)$, we have

$$P(A_n^\varepsilon) = P(A_n^\varepsilon | B_n^{\varepsilon_1}) P(B_n^{\varepsilon_1}) + P(A_n^\varepsilon | (B_n^{\varepsilon_1})^c) P((B_n^{\varepsilon_1})^c)$$

$$\leq P(B_n^{\varepsilon_1}) + P \left( C \{ \inf \mathbb{E} | Z_t^{+2} - Z_t^2 | \} \sum_{k=1}^N (k+1)(1 + \varepsilon_1)^k (1 - \delta)^k > \varepsilon \right)$$

$$\leq P(B_n^{\varepsilon_1}) + P(C_1 \inf \mathbb{E} | Z_t^{+2} - Z_t^2 | > \varepsilon) \to 0,$$

which demonstrates the convergence in probability of on the right-hand side of (59).

We now consider the second term $II$. Since $k > bN/(2p)$ and $\sup_u \sum_{j=1}^p a_j(u) \leq 1 - \delta$, it is straightforward to show

$$II \leq a_0(u_0) \sum_{k > bN/(2p)}^{\infty} (1 + (\mathbb{E} | Z_t^{+2} |)^k)(1 - \delta)^k$$

$$\leq a_0(u_0)(1 - \delta)^{bN/(2p)} \sum_{k=1}^{\infty} (1 + (\mathbb{E} | Z_t^{+2} |)^{k+bN/(2p)})(1 - \delta)^k.$$

Now it is straightforward to show that $II \xrightarrow{\mathcal{P}} 0$ with $b \to 0$, $bN \to \infty$ as $N \to \infty$. Altogether we obtain (56), and the desired result follows. $\square$

We note that the bounds given in Lemma A.7 are uniform for all $t_0 - bN/2 \leq t \leq t_0 + bN/2$, this is required to prove (52). As a byproduct of Lemma A.7, we have the following result.

COROLLARY 1. *Suppose the assumptions in Lemma* A.7 *hold. Then, for all* $t \in [t_0 - bN/2, t_0 + bN/2 - 1]$, *we have*

$$\mathbb{E}|\sigma_t^{+2}(u_0) - \tilde{\sigma}_t^{+2}(u_0)| \xrightarrow{\mathcal{P}} 0, \tag{60}$$

$$\inf \mathbb{E}|\tilde{\sigma}_t^{+2}(u_0) - \bar{\sigma}_t^2(u_0)| \xrightarrow{\mathcal{P}} 0, \tag{61}$$

*where* $b \to 0$, $bN \to \infty$ *as* $N \to \infty$, *and the expectations are defined in the same way as in Lemma* A.7.

PROOF. By using the expressions $\sigma_t^{+2}(u_0) = \bar{a}_{t_0,N}(0) + \sum_{j=1}^p \bar{a}_{t_0,N}(j) \times X_{t-j}^{+2}(u_0)$ and $\tilde{\sigma}_t^2(u) = a_0(u) + \sum_{j=1}^p a_j(u)\tilde{X}_{t-j}^2(u)$, and taking also into account that $\bar{\underline{a}}_{t_0,N} \xrightarrow{\mathcal{P}} \underline{a}(u_0)$, the desired result follows immediately from Lemma A.7. $\square$

In order to prove (52), we require the following inequalities.

Let us suppose $\sigma_x^2 = \alpha_0 + \sum_{j=1}^p \alpha_j x_j$, $\sigma_y^2 = \beta_0 + \sum_{j=1}^p \beta_j y_j$ with $\{\alpha_j\}$, $\{\beta_j\}$, $\{x_j\}$ and $\{y_j\}$ positive. Then, it can be shown that

$$
\begin{aligned}
&\left| \frac{z_x \sigma_x^2 x_i}{(\hat{\mu}_{t_0,N} + \sum_{j=1}^p x_j)^2} - \frac{z_x \sigma_y^2 y_i}{(\hat{\mu}_{t_0,N} + \sum_{j=1}^p y_j)^2} \right|^2 \\
&\leq \frac{K z_x^2 (A+B)^2}{\hat{\mu}_{t_0,N}} \left\{ 2\sum_{j=1}^p |x_j - y_j| + |\sigma_x^2 - \sigma_y^2| \right\},
\end{aligned}
\tag{62}
$$

where

$$A = \frac{\alpha_0}{\hat{\mu}_{t_0,N}} + \sum_{j=1}^p \alpha_j \quad \text{and} \quad B = \frac{\beta_0}{\hat{\mu}_{t_0,N}} + \sum_{j=1}^p \beta_j.$$

Similarly, we have

$$
\begin{aligned}
&\left| \frac{z_x \sigma_x^2 x_i}{(\hat{\mu}_{t_0,N} + \sum_{j=1}^p x_j)^2} - \frac{z_y \sigma_y^2 y_i}{(\hat{\mu}_{t_0,N} + \sum_{j=1}^p y_j)^2} \right|^2 \\
&\leq K(A+B)|z_x - z_y|^2 \\
&\quad + \frac{K z_x^2 (A+B)^2}{\hat{\mu}_{t_0,N}} \left\{ 2\sum_{j=1}^p |x_j - y_j| + |\sigma_x^2 - \sigma_y^2| \right\}.
\end{aligned}
\tag{63}
$$

We use these inequalities to prove the following result.

PROOF OF (52). By definition of Mallows metric, independence of the pairs $\{(Z_t^{+2}, Z_t^2)\}_t$, and that $\mathbb{E}(Z_t^{+2}) = 1$, we have

$$d_2\{\sqrt{bN}(\underline{r}_{t_0,N}^+ - \mathcal{R}_{t_0,N}^+ \bar{\underline{a}}_{t_0,N}), \sqrt{bN}(\tilde{\underline{r}}_N(u_0) - \tilde{\mathcal{R}}_N(u_0)\underline{a}(u_0))\}$$

$$
\begin{aligned}
(64) \quad &\leq (bN) \inf \mathbb{E}\{ (\underline{r}_{t_0,N}^+ - \mathcal{R}_{t_0,N}^+ \underline{\bar{a}}_{t_0,N}) - (\underline{\tilde{r}}_N(u_0) - \tilde{\mathcal{R}}_N(u_0)\underline{a}(u_0)) \}^2 \\
&\leq \frac{2}{bN} \sum_{j=1}^{p} \sum_{k=p}^{N} W\left(\frac{t_0-k}{b}\right)^2 \inf \mathbb{E}\bigg( \frac{(Z_k^{+2}-1)\sigma_k^{+2}(u_0)X_{k-i}^{+2}(u_0)}{[\hat{\mu}_{t_0,N} + \sum_{j=1}^{p} X_{k-j}^{+2}(u_0)]^2} \\
&\qquad\qquad - \frac{(Z_k^2-1)\tilde{\sigma}_k^2(u_0)\tilde{X}_{k-i}^2(u_0)}{[\hat{\mu}_{t_0,N} + \sum_{j=1}^{p} \tilde{X}_{k-j}^2(u_0)]^2} \bigg)^2,
\end{aligned}
$$

where the infimum is taken over all joint measures on $(Z_t^{+2}, Z_t^2)$. We now consider

$$
\begin{aligned}
\mathbb{E}\bigg( &\frac{(Z_k^{+2}-1)\sigma_k^{+2}(u_0)X_{t-i}^{+2}(u_0)}{[\hat{\mu}_{t_0,N} + \sum_{j=1}^{p} X_{k-j}^{+2}(u_0)]^2} - \frac{(Z_k^2-1)\tilde{\sigma}_k^2(u_0)\tilde{X}_{k-i}^2(u_0)}{[\hat{\mu}_{t_0,N} + \sum_{j=1}^{p} \tilde{X}_{k-j}^2(u_0)]^2} \bigg)^2 \\
&\leq 2(I+II),
\end{aligned}
$$

where

$$
I = \mathbb{E}\bigg( \frac{(Z_k^{+2}-1)\sigma_k^{+2}(u_0)X_{k-i}^{+2}(u_0)}{[\hat{\mu}_{t_0,N} + \sum_{j=1}^{p} X_{k-j}^{+2}(u_0)]^2} - \frac{(Z_k^{+2}-1)\tilde{\sigma}_k^{+2}(u_0)\tilde{X}_{k-i}^{+2}(u_0)}{[\hat{\mu}_{t_0,N} + \sum_{j=1}^{p} \tilde{X}_{k-j}^{+2}(u_0)]^2} \bigg)^2,
$$

$$
II = \mathbb{E}\bigg( \frac{(Z_k^{+2}-1)\sigma_k^{+2}(u_0)X_{k-i}^{+2}(u_0)}{[\hat{\mu}_{t_0,N} + \sum_{j=1}^{p} X_{k-j}^{+2}(u_0)]^2} - \frac{(Z_k^2-1)\tilde{\sigma}_k^2(u_0)\tilde{X}_{k-i}^2(u_0)}{[\hat{\mu}_{t_0,N} + \sum_{j=1}^{p} \tilde{X}_{k-j}^2(u_0)]^2} \bigg)^2.
$$

Studying first $I$, and using (63), we have

$$
\begin{aligned}
I &\leq K\mathbb{E}(Z_k^{+2}-1)^2 \\
&\quad \times \mathbb{E}\bigg( \frac{\sigma_k^{+2}(u_0)X_{k-i}^{+2}(u_0)}{[\hat{\mu}_{t_0,N} + \sum_{j=1}^{p} X_{k-j}^{+2}(u_0)]^2} - \frac{\tilde{\sigma}_k^{+2}(u_0)\tilde{X}_{k-i}^{+2}(u_0)}{[\hat{\mu}_{t_0,N} + \sum_{j=1}^{p} \tilde{X}_{k-j}^{+2}(u_0)]^2} \bigg)^2 \\
&\leq \frac{K\mathbb{E}(Z_k^{+2}-1)^2(A_1+B_1)^2}{\hat{\mu}_{t_0,N}} \\
&\quad \times \bigg\{ 2\sum_{j=1}^{p} \mathbb{E}|X_{k-j}^{+2}(u_0) - \tilde{X}_{k-j}^2(u_0)| + \mathbb{E}|\sigma_k^{+2}(u_0) - \tilde{\sigma}_k^{+2}(u_0)| \bigg\},
\end{aligned}
$$

where $A_1 = \frac{\hat{a}_{t_0,N}(0)}{\hat{\mu}_{t_0,N}} + \sum_{j=1}^{p} \hat{a}_{t_0,N}(j)$ and $B_1 = \frac{a_0(u_0)}{\hat{\mu}_{t_0,N}} + \sum_{j=1}^{p} a_j(u_0)$. Therefore, by using (54), (55) and (60), we have $I \xrightarrow{\mathcal{P}} 0$. Bounding $II$ by using (64), we have

$$
\begin{aligned}
II &\leq (A_1+B_1)\mathbb{E}|Z_k^{+2} - Z_k^2|^2 \\
&\quad + \frac{K\mathbb{E}(Z_k^{+2})(A_1+B_1)^2}{\hat{\mu}_{t_0,N}}
\end{aligned}
$$

$$\times \left\{ 2 \sum_{j=1}^{p} \mathbb{E}|\tilde{X}_{k-j}^{+2}(u_0) - \tilde{X}_{k-j}^{2}(u_0)| + \mathbb{E}|\tilde{\sigma}_k^{+2}(u_0) - \tilde{\sigma}_k^{2}(u_0)| \right\}.$$

Substituting the above bounds into (65), we have

$$d_2\{\sqrt{bN}(\underline{r}_{t_0,N}^+ - \mathcal{R}_{t_0,N}(u_0)^+ \underline{a}_{t_0,N}), \sqrt{bN}(\underline{\tilde{r}}_N(u_0) - \tilde{\mathcal{R}}_N(u_0)\underline{a}(u_0))\}$$
$$\leq \tilde{I} + \tilde{II},$$

where

$$\tilde{I} = \left[ \frac{4\mathbb{E}(Z_k^{+2} - 1)^2(A_1 + B_1)^2}{\hat{\mu}_{t_0,N}} \right.$$
$$\left. \times \left\{ 2 \sum_{j=1}^{p} \mathbb{E}|X_{k-j}^{+2}(u_0) - \tilde{X}_{k-j}^{+2}(u_0)| + \mathbb{E}|\sigma_k^{+2}(u_0) - \tilde{\sigma}_k^{+2}(u_0)| \right\} \right] \omega_N,$$

$$\tilde{II} = \frac{4\mathbb{E}(Z_k^{+2})(A_1 + B_1)^2}{\hat{\mu}_{t_0,N}}$$
$$\times \left\{ 2 \sum_{j=1}^{p} \inf \mathbb{E}|\tilde{X}_{k-j}^{+2}(u_0) - \tilde{X}_{k-j}^{2}(u_0)| + \inf \mathbb{E}|\tilde{\sigma}_k^{+2}(u_0) - \tilde{\sigma}_k^{2}(u_0)| \right\} \omega_N,$$

and $\omega_N = \frac{1}{bN} \sum_{j=1}^{p} \sum_{k=bN/2}^{bN/2} W(\frac{t_0 - k}{b})^2$. By using (54), (55) and (60), we have $\tilde{I} \xrightarrow{\mathcal{P}} 0$. By using (54), (56) and (61), we have $\tilde{II} \xrightarrow{\mathcal{P}} 0$. Altogether we obtain the required result. □

PROOF OF (53). We use the same methods as those in the proof of (52) to show that $d_2(\mathcal{R}_{t_0,N}^+, \tilde{\mathcal{R}}_N(u_0)) \xrightarrow{\mathcal{P}} 0$. Then, by using Lemma 8.3 in Bickel and Freedman (1981), and $\tilde{\mathcal{R}}_N(u_0) \xrightarrow{\mathcal{P}} \mathbb{E}[\mathcal{A}_t^{(\mu)}(u)]$ we have $\mathcal{R}_{t_0,N}^+ \xrightarrow{\mathcal{P}} \mathbb{E}[\mathcal{A}_t^{(\mu)}(u)]$, thus obtaining the desired result. □

We now have the necessary ingredients to prove Proposition 5.

PROOF OF PROPOSITION 5. We observe that

$$\sqrt{bN}(\underline{\hat{a}}_{t_0,N}^+ - \underline{\bar{a}}_{t_0,N})$$
$$= \sqrt{bN}(\mathcal{R}_{t_0,N}^+)^{-1}(\underline{r}_{t_0,N}^+ - \mathcal{R}_{t_0,N}^+ \underline{\bar{a}}_{t_0,N}).$$

Now, since by (53) we have $\mathcal{R}_{t_0,N}^+ \xrightarrow{\mathcal{P}} \mathbb{E}[A_t^{(\mu)}(u_0)]$, we can replace in the above $\mathcal{R}_{t_0,N}^+$ with $\mathbb{E}[A_t^{(\mu)}(u_0)]$, and then use the delta method and (52) to get the required result. □

**Acknowledgments.** Piotr Fryzlewicz would like to thank David Hand and Theofanis Sapatinas for financial support while visiting Nicosia to carry out part of this work. We also thank Esther Ruiz for interesting discussions. Most of the research was conducted while Suhasini Subba Rao was at the Universität Heidelberg. We are grateful to the co-editor (Jianqing Fan), an Associate Editor, and the two reviewers whose valuable comments and suggestions led to a significant improvement of this paper.

## REFERENCES

BERA, A. K. and HIGGINS, M. L. (1993). ARCH models: Properties, estimation and testing. *J. Econom. Surv.* **7** 305–366.

BHATTACHARYA, R. N., GUPTA, V. K. and WAYMIRE, E. (1983). The Hurst effect under trend. *J. Appl. Probab.* **20** 649–662. MR0713513

BICKEL, P. J. and FREEDMAN, D. A. (1981). Some asymptotic theory for the bootstrap. *Ann. Statist.* **9** 1196–1217. MR0630103

BOLLERSLEV, T. (1986). Generalized autoregressive conditional heteroskedasticity. *J. Econometrics* **31** 307–327. MR0853051

BOSE, A. and MUKHERJEE, K. (2003). Estimating the ARCH parameters by solving linear equations. *J. Time Ser. Anal.* **24** 127–136. MR1965808

DAHLHAUS, R. (1997). Fitting time series models to nonstationary processes. *Ann. Statist.* **25** 1–37. MR1429916

DAHLHAUS, R. and SUBBA RAO, S. (2006). Statistical inference for time-varying ARCH processes. *Ann. Statist.* **34** 1075–1114. MR2278352

ENGLE, R. F. (1982). Autoregressive conditional heteroscedasticity with estimates of the variance of United Kingdom inflation. *Econometrica* **50** 987–1008. MR0666121

FAN, J. and YAO, Q. (2003). *Nonlinear Time Series.* Springer, New York. MR1964455

FAN, J., JIANG, J., ZHANG, C. and ZHOU, Z. (2003). Time-dependent diffusion models for term structure dynamics. *Statist. Sinica* **13** 965–992. MR2026058

FRANKE, J. and KREISS, J.-P. (1992). Bootstrapping stationary autoregressive moving average models *J. Time Ser. Anal.* **13** 297–317. MR1173561

FRYZLEWICZ, P., SAPATINAS, T. and SUBBA RAO, S. (2006). A Haar–Fisz technique for locally stationary volatility estimation. *Biometrika* **93** 687–704. MR2281147

GIRAITIS, L., KOKOSZKA, P. and LEIPUS, R. (2000). Stationary ARCH models: Dependence structure and central limit theorem. *Econometric Theory* **16** 3–22. MR1749017

GIRAITIS, L., LEIPUS, R. and SURGAILIS, D. (2005). Recent advances in ARCH modelling. In *Long Memory in Economics* (A. Kirman and G. Teyssiere, eds.) 3–38. Springer, Berlin. MR2265054

GIRAITIS, L. and ROBINSON, P. (2001). Whittle estimation of GARCH models. *Econometric Theory* **17** 608–631. MR1841822

HALL, P. and HEYDE, C. (1980). *Martingale Limit Theory and its Applications.* Academic Press, New York. MR0624435

HART, J. (1996). Some automated methods for smoothing time-dependent data. *J. Nonparametr. Statist.* **6** 115–142. MR1383047

HORVÁTH, L., KOKOSZKA, P. and TEYSSIÉRE, G. (2001). Empirical process of the squared residuals of an ARCH process. *Ann. Statist.* **29** 445–469. MR1863965

HORVÁTH, L. and LIESE, F. (2004). $L_p$-estimators in ARCH models. *J. Statist. Plann. Inference* **119** 277–309. MR2019642

LING, S. (2007). Self-weighted and local quasi-maximum likelihood estimators for ARMA-GARCH/IGARCH models. *J. Econom.* **140** 849–873.

MERCURIO, D. and SPOKOINY, V. (2004a). Statistical inference for time-inhomogeneous volatility models. *Ann. Statist.* **32** 577–602. MR2060170

MERCURIO, D. and SPOKOINY, V. (2004b). Estimation of time dependent volatility via local change point analysis. Preprint.

MIKOSCH, T. and STĂRICĂ, C. (2000). Is it really long memory we see in financial returns? In *Extremes and Integrated Risk Management* (P. Embrechts, ed.) 149–168. Risk Books, London.

MIKOSCH, T. and STĂRICĂ, C. (2003). Long-range dependence effects and ARCH modelling. In *Theory and Applications of Long Range Dependence* (P. Doukhan, G. Oppenheim and M. S. Taqqu, eds.) 439–459. Birkhäuser, Boston. MR1957503

MIKOSCH, T. and STĂRICĂ, C. (2004). Non-stationarities in financial time series, the long-range dependence, and the IGARCH effects. *Rev. Econ. Statist.* **86** 378–390.

PAPARODITIS, E. and POLITIS, D. N. (2007). Resampling and subsampling for financial time series. In *Handbook of Financial Time Series* (T. Andersen, R. A. Davis, J.-P. Kreiss and T. Mikosch, eds.). Springer, New York. To appear.

SHEPHARD, N. (1996). Statistical aspects of ARCH and stochastic volatility. In *Time Series Models in Econometric, Finance and Other Fields* (D. R. Cox, D. V. Hinkley and O. E. Barndorff-Nielsen, eds.) 1–67. Chapman and Hall, London.

STRAUMANN, D. (2005). *Estimation in Conditionally Heteroscedastic Time Series Models*. Springer, New York. MR2142271

STĂRICĂ, C. (2003). Is GARCH(1, 1) as good a model as the Nobel prize accolades would imply? Preprint.

STĂRICĂ, C. and GRANGER, C. W. J. (2005). Non-stationarities in stock returns. *Rev. Econ. Statist.* **87** 503–522.

SUBBA RAO, S. (2006). On some nonstationary, nonlinear random processes and their stationary approximations. *Adv. in Appl. Probab.* **38** 1153–1172. MR2285698

TAYLOR, S. C. (1986). *Modelling Financial Time Series*. Wiley, Chichester.

P. FRYZLEWICZ
SCHOOL OF MATHEMATICS
UNIVERSITY OF BRISTOL
BRISTOL BS8 1TW
UNITED KINGDOM
E-MAIL: p.z.fryzlewicz@bristol.ac.uk

T. SAPATINAS
DEPARTMENT OF MATHEMATICS
  AND STATISTICS
UNIVERSITY OF CYPRUS
CY 1678 NICOSIA
CYPRUS
E-MAIL: t.sapatinas@ucy.ac.cy

S. SUBBA RAO
DEPARTMENT OF STATISTICS
TEXAS A&M UNIVERSITY
COLLEGE STATION, TEXAS 77843
USA
E-MAIL: suhasini.subbarao@stat.tamu.edu